\documentclass[reqno,11pt]{amsart}
\usepackage{amsmath,amsfonts,amssymb,amsthm}

\usepackage{latexsym,bm,graphicx}
\usepackage{mathrsfs}
\usepackage{color}

\newtheorem{thm}{Theorem}[section]
\newtheorem{prop}[thm]{Proposition}
\newtheorem{lem}[thm]{Lemma}
\newtheorem{cor}[thm]{Corollary}

\theoremstyle{definition}
\theoremstyle{remark}

\newtheorem{defn}[thm]{Definition}
\newtheorem{remark}[thm]{Remark}

\numberwithin{equation}{section}

\newcommand{\ls}{\leqslant}
\newcommand{\gs}{\geqslant}

\newcommand{\rd}{\mathrm{d}}
\newcommand{\rvol}{\mathrm{vol}}
\newcommand{\bR}{\mathbb{R}}
\newcommand{\eps}{\epsilon}
\newcommand{\rLip}{\mathrm{Lip}}
\newcommand{\cH}{\mathcal{H}}
\newcommand{\cL}{\mathcal{L}}
\newcommand{\SL}{\mathscr{L}}
\newcommand{\tM}{\widetilde{M}}
\newcommand{\Reg}{\mathrm{Reg}}
\newcommand{\tih}{\tilde{h}}
\newcommand{\tf}{\tilde{f}}
\newcommand{\tg}{\tilde{g}}
\newcommand{\bH}{\mathbb{H}}

\def\Xint#1{\mathchoice
{\XXint\displaystyle\textstyle{#1}}%
{\XXint\textstyle\scriptstyle{#1}}%
{\XXint\scriptstyle\scriptscriptstyle{#1}}%
{\XXint\scriptscriptstyle\scriptscriptstyle{#1}}%
\!\int}
\def\XXint#1#2#3{{\setbox0=\hbox{$#1{#2#3}{\int}$ }
\vcenter{\hbox{$#2#3$ }}\kern-.6\wd0}}

\def\dashint{\Xint-}

\title{maximal bottom of spectrum or volume entropy rigidity in Alexandrov geometry}
\author{Yin Jiang}
\address{Mathematics Department, Capital Normal University, Beijing P.R.C.}
\email[Yin Jiang]{jiangyinwd@163.com}

\begin{document}
\maketitle

\bibliographystyle{amsplain}

\begin{abstract}
In \cite{LiWang2001complete1,LiWang2001complete2}, Li-Wang proved a splitting theorem for an n-dimensional Riemannian manifold with $Ric\gs -(n-1)$ and the bottom of spectrum $\lambda_0(M)=\frac{(n-1)^2}{4}$. For an n-dimensional compact manifold $M$ with $Ric\gs -(n-1)$ with the volume entropy $h(M)=n-1$, Ledrappier-Wang \cite{LeW2010volent} proved that the universal cover $\tM$ is isometric to the hyperbolic space $\bH^n$. We will prove analogue theorems for Alexandrov spaces.

\end{abstract}
\maketitle
\section{Introduction}
Let $(M,|\cdot \cdot|)$ be an n-dimensional Alexandrov space with curvature $\gs -1$, $\partial M=\emptyset$. Denote $\rvol$ its n-dimensional Hausdorff measure. The goal of this paper is to establish two rigidity theorems on Alexandrov spaces.

Given a Lipschitz function $f: M \to \bR$, the pointwise Lipschitz constant of $f$ at $x$ is defined by
\[
\rLip f(x):=\limsup_{y \to x} \frac{|f(y)-f(x)|}{|xy|}.
\]
Denote by $\rLip_c(M)$ the set of Lipschitz functions with compact support in $M$. Suppose $M$ is non-compact, the bottom of the $L^2$-spectrum of the Laplacian on $M$ can be characterized by
\[
\lambda_0(M)=\inf_{f\in \rLip_c(M)} \frac{\int_M (\rLip f)^2 \rd \rvol}{\int_M f^2 \rd \rvol}.
\]
It's well known that (see e.g. Theorem 5 of \cite{sturm1994analysis}) that
\begin{equation}\label{ils:eigenvol}
\sqrt{\lambda_0(M)} \ls \frac12 \limsup_{R\to \infty} \frac{\ln \rvol(B_p(R))}{R}.
\end{equation}
By the Bishop volume comparison, we have
\begin{equation}\label{ils:lambda0}
\lambda_0(M) \ls \frac{(n-1)^2}{4}.
\end{equation}

When $(M,g)$ is a smooth Riemannian manifold, Li-Wang \cite{LiWang2001complete1,LiWang2001complete2} proved the following theorems
\begin{thm}[Theorem 0.5 and 0.6, \cite{LiWang2001complete2}]\label{thm:ngs3}
Let $M^n$ be a complete n-dimensional manifold. Suppose that
\begin{equation}\label{igs:Ric}
Ric_M \gs -(n-1)
\end{equation}
and
\begin{equation}\label{igs:eig}
\lambda_0(M) = \frac{(n-1)^2}{4}.
\end{equation}
If $n\gs 4$, then either:

(1) $M$ has only one end; or

(2) $M=\bR\times N$ with the warped product metric
\begin{equation}
\rd s_M^2=\rd t^2+\exp(2t) \rd s_N^2,
\end{equation}
where $N$ is a compact manifold with non-negative Ricci curvature.

If $n=3$, besides (1) and (2), we have another case:

(3) $M=\bR \times N$ with the warped product metric
\begin{equation}
\rd s_M^2=\rd t^2+\cosh^2 t \rd s_N^2,
\end{equation}
where $N^2$ is a compact manifold with its Gaussian curvature bounded below by $K_N \gs -1$.
\end{thm}

For $n=2$, there is no splitting theorem for two infinite volume ends, see section 3 of \cite{LiWang2001complete2}. However, we still have the following:
\begin{thm}[Theorem 0.7, \cite{LiWang2001complete2}]\label{thm:n=2}
Let $M^2$ be a complete 2-dimensional manifold. Suppose that $K_M \gs -1$ and $\lambda_0(M)=\frac14$, then either:

(1) $M$ has no finite volume end; or

(2) $M=\bR \times S^1$ with the warped product metric
\begin{equation}
\rd s_M^2=\rd t^2+\exp(2t) \rd \theta^2.
\end{equation}
\end{thm}

Let $M$ be a compact Riemannian manifold, denote by $\tM$ its universal cover. If $Ric_M \gs -(n-1)$ and $\lambda_0(\tM)=\frac{(n-1)^2}{4}$, from the above theorems, we know that $\tM$ has at most one end. Wang \cite{Wang2008harmonic} proved that $\tM$ must be isometric to the hyperbolic space $\bH^n$. Later, the condition $\lambda_0(\tM)=\frac{(n-1)^2}{4}$ was weakened by Ledrappier-Wang \cite{LeW2010volent}. The volume entropy of a compact manifold $M$ is defined by
\begin{equation}
h(M)=\lim_{R \to \infty} \frac{\ln \rvol(B_R(\tilde{p}))}{R}, \quad \tilde{p} \in \tM.
\end{equation}
(for the existence of the limit, see \cite{M1979ventropy}).
By Bishop volume comparison, for any compact $n$-dimensional manifold $M$ with $Ric_M \gs -(n-1)$, $h(M) \ls n-1$.
\begin{thm}[\cite{LeW2010volent}]\label{thm:LeW}
Let $M$ be a compact, $n$-dimensional Riemannian manifold with $Ric \gs -(n-1)$. If $h(M)=n-1$, then $\tM$ is isometric to the hyperbolic space $\bH^n$.
\end{thm}
\begin{remark}
This is a generalization of Wang's theorem. Since if $\lambda_0(M)=\frac{(n-1)^2}{4}$, by (\ref{ils:eigenvol}), we have $h(M)=(n-1)$. Recently, Chen-Rong-Xu \cite{CRX2016Quati} have proved a quantitative version for Theorem \ref{thm:LeW}.
\end{remark}

If $M$ is a compact n-dimensional Alexandrov space with curvature $\gs -1$, by Bishop volume comparison, the volume entropy $h(M)\ls (n-1)$. In view of the theorems above, do we have any rigidity for non-compact Alexandrov spaces with curvature $\gs -1$ satisfying $\lambda_0=\frac{(n-1)^2}{4}$? Or compact Alexandrov spaces with curvature $\gs -1$ satisfying $h(M)=(n-1)$? In this paper, we will prove analogue theorems for Alexandrov spaces.
\begin{thm}\label{thm:eigensplitting}
Suppose $n\gs 4$. Let $M$ be a non-compact, n dimensional Alexandrov space with curvature $\gs -1$, $\partial M=\emptyset$. If $\lambda_0(M)=\frac{(n-1)^2}{4}$, then either

(1) $M$ has only one end; or

(2) $M$ splits as $M=\bR \times_{e^t} N$, where $N$ is a compact Alexandrov space with non-negative curvature.
\end{thm}
For $n=2,3$, we have the following theorem
\begin{thm}\label{thm:An23}
Suppose $n=2 \text{ or } 3$. Let $M$ be a non-compact, n-dimensional Alexandrov space with curvature $\gs -1$, $\partial M=\emptyset$. If $\lambda_0(M)=\frac{(n-1)^2}{4}$, then either

(1) $M$ has no finite volume end; or

(2) $M$ splits as $M=\bR \times_{e^t} N$, where $N$ is a compact Alexandrov space with non-negative curvature.
\end{thm}
Note that when $n=3$, our theorem is weaker than Theorem \ref{thm:ngs3}. Since if $M$ has at least two infinite volume ends , we don't know whether $M$ is a warped product like case (3) in Theorem \ref{thm:ngs3}, see Remark \ref{rem:3weaker}.

We will also prove a version of Theorem \ref{thm:LeW} for Alexandrov spaces.
\begin{thm}\label{thm:volent}
Let $M$ be a compact, $n$-dimensional Alexandrov space with curvature $\gs -1$. If the volume entropy $h(M)=n-1$, then $M$ is a hyperbolic manifold.
\end{thm}

As an immediate corollary of Theorem \ref{thm:volent}, we have
\begin{cor}
Let $M$ be a compact, $n$-dimensional Alexandrov space with curvature $\gs -1$. If $\lambda_0(\tM)=\frac{(n-1)^2}{4}$, then $M$ is a hyperbolic manifold.
\end{cor}

Before describing our approach, let us recall the proof of Theorem \ref{thm:ngs3}.
Suppose $M$ has at least two ends. If $M$ has at least two infinite volume ends and $\lambda_0(M)>0$, then we can construct a non-constant, bounded harmonic function on $M$. If in addition $Ric_M\gs -(n-1)$ and $\lambda_0(M)\gs n-2$, by Bochner formula and decay estimates for harmonic functions, Li-Wang \cite{LiWang2001complete1} proved that $M$ must splits as case (3) in Theorem \ref{thm:ngs3}, then $\lambda_0(M)=n-2$.

If $n\gs 4$, then $\frac{(n-1)^2}{4} >n-2$, $M$ has at least one finite volume end $E$. Ji-Li-Wang \cite{JLW2009ends} proved that the Busemann function with respect to the ray to the infinity of $E$ satisfies $\Delta b=n-1$. Then $b$ is smooth and $|\nabla b|=1$, $b$ has no critical point. It follows that $M$ is homeomorphic to $R \times N$ for some manifold $N$. By Bochner formula, they get the explicit form of the Hessian of $b$ and proved that $M=\bR \times_{e^t} N$. The proof of Theorem \ref{thm:n=2} is just the same.

For Theorem \ref{thm:LeW}, Liu \cite{Liu2011short} constructed a Busemann function on $\tM$ such that $\Delta b=n-1$ and $|\nabla b|=1$. By the argument as above $\tM=\bR \times_{e^t} N$. Since the sectional curvature of $M$ is bounded, by a theorem of \cite{Wang2008harmonic}, $\tM$ is isometric to $\bH^n$.

Our proof of Theorem \ref{thm:eigensplitting}, \ref{thm:An23} and \ref{thm:volent} are basically along the Line of the argument above. However, for Alexandrov spaces, due to the lack of smoothness of the boundaries of ends, harmonic functions are not necessarily continuous up to the boundary. We should rely on the theory of Dirichlet problem on metric spaces with a doubling measure and satisfying $(1,p)$ Poincar\'e inequality for $p>1$. Following the approach by harmonic functions developed by Li-Wang, we can prove that if $\lambda_0(M)>n-2$, then $M$ doesn't have two infinite volume ends.

For $n\gs 4$, since $\frac{(n-1)^2}{4} >n-2$, then $M$ has at least one end with finite volume. Following Ji-Li-Wang's proof, we can get a semiconcave function $b:M \to \bR$ such that $|\nabla x b|=1$ for a.e. $x\in M$ and $\SL_b=(n-1)\cdot \rvol$ (see section 2 for the definition of the Laplacian $\SL$). Under the condition of Theorem \ref{thm:volent}, following Liu's proof, we can also get such a function on $\tM$.

Similar to the non-negative curvature case, Alexander-Bishop \cite{AleBi05cone} proved that the existence of an affine function is equivalent to the splitting of an Alexandrov space. For our purpose, we just mention a particular case of this theorem.
\begin{defn}
We say $f:M \to \bR$ is a $-1$-affine function, if for any unit speed geodesic $\gamma(t)$,
\begin{equation}
[f\circ \gamma(t)]''-f\circ \gamma(t)=0.
\end{equation}
\end{defn}

\begin{thm}[\cite{AleBi05cone}]\label{thm:conesplit}
Let $M$ be an n-dimensional Alexandrov space with curvature $\gs -1$, $\partial M=\emptyset$. If $M$ carries a non-constant $-1$-affine function $f:M \to \bR$. Then $|\nabla_x f|^2-f^2$ is a constant and $M$ is a warped product. If in addition $|\nabla_x f|^2-f^2=0$, then $M=\bR \times_{e^t} N$, where $N$ is an Alexandrov space with non-negative curvature.
\end{thm}

If $M$ is a manifold, $e^b$ is just the $-1$ affine function.
For Alexandrov spaces, due to the lack of regularity for functions with constant Laplacian, it's not easy to see that $b$ has no critical points. For theorem \ref{thm:eigensplitting}, let $x\in M$, one may consider two asymptotic rays from $x$ with respect to a line on $M$. However, in general these two rays on the warped product $\bR \times_{e^t} N$ don't form a line. By studying the gradient flow, we will prove the following general splitting theorem.
\begin{thm}\label{thm:general}
Let $M$ be a non-compact, n-dimensional Alexandrov space with curvature $\gs -1$, $\partial M=\emptyset$. If there exists a semiconcave function $b:M \mapsto \bR$ satisfying:

\begin{equation}\label{ieq:ae0}
(i) \qquad |\nabla_x b|=1 \text{ for } \cH^n-a.e. \quad x\in M.
\end{equation}

\begin{equation}\label{ieq:Lb}
(ii) \qquad \SL_b=(n-1)\cdot \rvol.
\end{equation}
Then $f=e^b$ is $-1$-affine and $M$ splits as $M=\bR \times_{e^t} N$, where $N$ is an $n-1$ dimensional Alexandrov space with non-negative curvature.
\end{thm}
For Theorem \ref{thm:volent}, by an argument of Chen-Rong-Xu \cite{CRX2016Quati}, the warped product (i.e. $\tM$) must be $\bH^n$.
We will also discuss the obstacles to generalize our argument to $RCD^*(K,N)$ spaces, see Remark \ref{rem:RCD}.

The rest of the paper is organized as follows. In section 2, we recall some necessary materials for Alexandrov spaces, including gradient flow, theory of Dirichlet problem. We will also prove a refined version of localized Bochner inequality. In section 3, we will prove Theorem \ref{thm:general}. In section 4, we will prove theorem \ref{thm:eigensplitting} and  \ref{thm:An23}. In section 5, we will prove Theorem \ref{thm:volent}.

{\noindent\bf{Acknowledgements}.} We are grateful to Xiantao Huang, RenJin Jiang, Shicheng Xu and Huichun Zhang for helpful discussions. We also thank Xiaochun Rong for helpful suggestions.

\section{Preliminaries}
\subsection{Preliminaries on Alexandrov spaces}
In this section, we review the definition of Alexandrov spaces with curvature bounded below and some properties. These definitions and results are mainly taken from \cite{burago1992ad}, \cite{otsu1994riemannian} and \cite{burago2001course}.

Let $(M,|\cdot \cdot|)$ be a metric space. A rectifiable curve $\gamma$ connecting two points $p,q$ is called a geodesic if its length is equal to $|pq|$ and it has unit speed. A metric space is called a geodesic space if any two points $p,q \in M$ can be connected by a geodesic.
Denote by $M_k^2$ the simply connected 2-dimensional space form of constant curvature $k$. Given three points $p,q,r$ in a geodesic space $M$, we can take a comparison triangle $\Delta \tilde{p}\tilde{q}\tilde{r}$ in $M^2_k$, such that
\[
d(\tilde{p},\tilde{q})=|pq|,d(\tilde{p},\tilde{r})=|pr|,d(\tilde{q},\tilde{r})=|qr|.
\]
If $k>0$, we add the assumption $|pq|+|pr|+|qr|<2\pi/\sqrt{k}$. The angle $ \widetilde{\angle}_k pqr: =\angle \tilde{p}\tilde{q}\tilde{r}$ is called comparison angle.

\begin{defn}
A geodesic space $M$ is called an Alexandrov space with curvature $\gs k$ if it's locally compact and for any point $x\in M$, there exists a neighborhood $U_x$ such that, for any four different points $p,a,b,c$ in $U_x$, we have
\[
\widetilde{\angle}_k abp +\widetilde{\angle}_k bpc +\widetilde{\angle}_k cpa \ls 2\pi.
\]
\end{defn}
The Hausdorff dimension of an Alexandrov space is always an integer. Let $M$ be an n-dimensional Alexandrov space with curvature $\gs k$. Denote by $\cH^n$ the n-dimensional Hausdorff measure.
Given any two geodesics $\gamma(t)$ and $\eta(s)$ with $\gamma(0)=\eta(0)=p$, the angle
\[
\angle(\gamma^+(0),\eta^+(0)): =\lim_{t,s\to 0} \widetilde{\angle}_k \gamma(t)p\eta(s)
\]
is well defined.

We say $\eta(t)$ is equivalent to $\gamma(t)$ if $\angle(\gamma^+(0),\eta^+(0))=0$, denote by $\Sigma'_p$ the set of equivalent classes of geodesic $\gamma(t)$ with $\gamma(0)=p$. The space of directions $\Sigma_p$ is the completion of metric space $(\Sigma'_p,\angle)$.

The tangent cone at $p$, $T_p$, is the Euclidean cone over $\Sigma_p$, it's an Alexandrov space with curvature $\gs 0$. For any two vectors $u,v\in T_p$. The "scalar product" (see section 1 of \cite{petrunin2007semiconcave}) is defined by
\[
\langle u,v \rangle=|u||v|\cos \angle(u,v).
\]
The distance $|uv|$ is defined by the law of cosines
\begin{equation}\label{eq:lawcos}
|uv|^2=|u|^2+|v|^2-2|u||v|\cos \angle(u,v).
\end{equation}

For each point $x\neq p$, we denote by $\Uparrow_p^x$ the set of directions at $p$ corresponding to all geodesics connecting $p$ to $x$. The symbol $\uparrow_p^x$ denotes the direction at $p$ corresponding to some geodesic $px$.
Given a direction $\xi \in \Sigma_p$, it's possible that there exists no geodesic $\gamma(t)$ starting at $p$ with $\gamma^+(0)=\xi$. However, it's shown in \cite{perelmanquasi} that for $p\in M$ and any direction $\xi \in \Sigma_p$, there exists a quasi-geodesic $\gamma: [0,+\infty) \to M$ with $\gamma(0)=p$ and $\gamma^+(0)=\xi$.

The exponential map $\exp_p: T_p \to M$ is defined by Petrunin \cite{petrunin1998parallel} as follows. $\exp_p(o_p)=p$ and for any $v\in T_p \backslash \{o_p\}$, $\exp_p(v)$ is a point on some quasi-geodesic of length $|v|$ starting from $p$ along direction $\frac{v}{|v|}\in \Sigma_p$. If the quasi-geodesic is not unique, we fix some one of them as the definition of $\exp_p(v)$.

A point $p$ in an n-dimensional Alexandrov space $M$ is said to be regular if its tangent cone $T_p$ is isometric to $\bR^n$ with standard metric. Denote by $\Reg(M)$ the set of regular points.

\begin{defn}
We say that a function $u$ is differentiable at $x\in \text{Reg(M)}$, if there exists a vector in $T_x$, denoted by $\nabla u(x)$, such that for any geodesic $\gamma(t)$ with $\gamma(0)=x$,
\[
u(\gamma(t))=u(x)+\langle \nabla u(x),\gamma^+(0)\rangle t +o(t).
\]
\end{defn}
The Rademacher theorem, in the framework of metric measure space with a doubling measure and a Poincar\'e inequality for upper gradient, was proved by Cheeger \cite{cheeger1999diff}. In \cite{bertrand2008existence}, Bertrand proved it in Alexandrov spaces via a simple argument. It says that a locally Lipschitz function $u$ is differentiable almost everywhere with respect to $\cH^n$ in $M$.

\subsection{Semiconcave functions and gradient curves}
Next, we introduce $\lambda$-concave functions and semi-concave functions. These definitions and results are mainly taken from section 1 and 2 of \cite{petrunin2007semiconcave}.
\begin{defn}
Let $M$ be an n-dimensional Alexandrov space without boundary and $U \subset M$ be an open subset. A locally Lipschitz function $f: U \mapsto \bR$ is called $\lambda$-concave if for any geodesic $\gamma(t)$ in $U$, the function $f\circ \gamma(t)-\lambda t^2/2$ is concave.
\end{defn}
A function $f:M \mapsto \bR$ is called semiconcave if for any point $x \in M$, there is a neighborhood $U_x \ni x$ and $\lambda \in \bR$ such that the restriction $f|_{U_x}$ is $\lambda$-concave. Given a semiconcave function $f:M \mapsto \bR$, its differential $d_pf$ is well defined for each point $p\in M$.
Let $\varphi:\bR \to \bR$ be a continuous function. A function $f:M \mapsto \bR$ is called $\varphi(f)$-concave if for any point $x\in M$ and $\eps>0$, there is a neighborhood $U_x \ni x$ such that $f|_{U_x}$ is $(\varphi \circ f(x)+\eps)$-concave.

Note that any semiconcave function is locally Lipschitz. The gradient vector at any point $x$, $\nabla_x f$ is well defined. If $d_pf(v)\ls 0$ for all $v\in T_p$, then $\nabla_p f=o_p$; Otherwise,
\[
\nabla_p f=d_pf(\xi_{\max})\cdot \xi_{\max},
\]
where $\xi_{\max} \in \Sigma_p$ is the (necessarily unique) unit vector for which $d_p f$ attains its maximum.

Denote by $\rLip_{loc}(\Omega)$ the set of locally Lipschitz continuous functions on $\Omega$. Let $u\in Lip_{loc}(\Omega)$, the pointwise Lipschitz constant of $u$ at $x$ are defined by
\[
\rLip f(x):=\limsup_{y \to x} \frac{|f(y)-f(x)|}{|xy|}.
\]
For the gradient, we have the following proposition:
\begin{lem}[Proposition 2.4, \cite{zhang2012yau}]\label{lem:graequal}
Let $f: U \to \bR$ be a semiconcave function. If $f$ is differentiable at $x$, then we have
\begin{equation}
|\nabla f(x)|=|\nabla_x f|=\rLip u(x).
\end{equation}
\end{lem}

Next we introduce the gradient curves of semiconcave functions.
\begin{defn}
Let $f:M \to \bR$ be a semiconcave function, a curve $\alpha(t)$ is called $f$-gradient curve if for any $t$,
\[
\alpha^+(t)=\nabla_{\alpha(t)} f.
\]
\end{defn}
The next proposition states the existence and uniqueness of gradient curves.
\begin{prop}[Propostion 2.1.2, \cite{petrunin2007semiconcave}]
Given a semiconcave function $f:M \to \bR$ and a point $p\in M$, there is a unique gradient curve $\alpha:[0,\infty) \to M$ such that $\alpha(0)=p$.
\end{prop}
A limit of gradient curves is a gradient curve for the limit function, i.e.
\begin{prop}[Proposition 2.1.5, \cite{petrunin2007semiconcave}]\label{prop:gradcurv}
Let $p_n \to p$, let $\alpha_n:[0,\infty) \to M$ be the sequence of f-gradient curves with $\alpha_n(0)=p_n$ and let $\alpha:[0,\infty)\to M$ be the $f$-gradient curve with $\alpha(0)=p$. Then $\alpha_n \to \alpha$ as $n \to \infty$.
\end{prop}
Next we introduce the gradient flow.
\begin{defn}
Let $f:M \to \bR$ be a semiconcave function. We define the $f$-gradient flow to be the one parameter family of maps
\[
\Phi_f^t: M \to M, \quad \Phi_f^t(p)=\alpha_p(t),
\]
where $t\gs 0$ and $\alpha_p:[0,\infty)\to M$ is the $f$-gradient curve which starts at $p$.
\end{defn}

\subsection{Sobolev spaces and measure valued Laplacian}
Let $\Omega$ be a domain in $M$, the Sobolev spaces $W^{1,2}(\Omega)$ is well defined (see, for example \cite{KMS2001sobolev}).
For a locally Lipschitz function $u$, its $W^{1,2}(\Omega)$-norm is defined by
\[
\|u\|_{W^{1,2}(\Omega)}:=\|u\|_{L^2(\Omega)}+\|\rLip u\|_{L^2(\Omega)}.
\]
Sobolev spaces $W^{1,2}(\Omega)$ is defined by the closure of the set
\[
\{u\in \rLip_{loc}(\Omega)| \|u\|_{W^{1,2}(\Omega)}<\infty\}
\]
under the $W^{1,2}(\Omega)$-norm. Denote by $\rLip_c(\Omega)$ the set of Lipschitz functions with compact support in $\Omega$. $W^{1,2}_0(\Omega)$ is defined by the closure of $\rLip_c(\Omega)$ under the $W^{1,2}(\Omega)$-norm. This coincides with the definitions in \cite{cheeger1999diff}. We say $u\in W^{1,2}_{loc}(\Omega)$ if $u\in W^{1,2}(\Omega')$ for any bounded, open subset $\Omega' \Subset \Omega$. According to \cite{KMS2001sobolev} (see also Theorem 4.47 of \cite{cheeger1999diff}), the "derivative" $\nabla u$ is well-defined for all $u\in W^{1,2}(\Omega)$. $W^{1,2}(\Omega)$ is reflexive according to Theorem 4.48 of \cite{cheeger1999diff}.

Given a function $u \in W^{1,2}_{loc}(\Omega)$, a functional $\SL_u$ is defined on $\rLip_c(\Omega)$ by
\[
\SL_u(\phi)=-\int_{\Omega} \langle \nabla u, \nabla \phi \rangle \rd \rvol, \quad \forall \phi \in \rLip_c(\Omega).
\]

By a standard argument, we can prove the following Lemma:
\begin{lem}\label{plem:conver}
Let $u_n \in W^{1,2}_{loc}(\Omega)$ and $u\in L^2_{loc}(\Omega)$. If for any bounded, open subset $\Omega' \Subset \Omega$, there exists a constant $C(\Omega')$ such that $\|u_n\|_{W^{1,2}(\Omega')} \ls C(\Omega')$ and $u_n$ converge to $u$ strongly in $L^2(\Omega')$, then $u\in W^{1,2}_{loc}(\Omega)$ and for any $\phi \in Lip_c(\Omega)$,
\[
\SL_{u_n}(\phi) \to \SL_u (\phi) \text{ as } n \to \infty.
\]
\end{lem}

Let $f\in L^2_{loc}(\Omega)$, if for any non-negative $\phi \in \rLip_c(\Omega)$,
\[
\SL_u(\phi)\ls \int_{\Omega} f\phi \rd \rvol,
\]
then we say $\SL_u\ls f\cdot \rvol$. In this case, according to \cite{Ho1989pde}, $\SL_u$ is a signed Radon measure. Denote its Lebesgue decompostion by
\[
\SL_u=\Delta^{ac} u \cdot \rvol +\Delta^s u,
\]
where $\Delta^{ac} u$ is the density of the absolutely continuous part and $\Delta^s u$ is the singular part. We have that
\[
\Delta^{ac} u(x) \ls f(x) \text{ for } \cH^n \  a.e. \  x \in \Omega \text{ and } \Delta^s u \ls 0.
\]
For a semiconcave function $f:M \to \bR$, it was proved by Perelman \cite{pereDC} that for a.e. $p \in \Reg(M)$, there exists a quadratic form $H_p f$ on $T_x$ such that for any geodesic $\gamma(t)$ with $\gamma(0)=p$, we have
\begin{equation}\label{peq:perhess}
f\circ \gamma(t)-f(p)=d_p f(\gamma'(0))t+ \frac12 H_p f(\gamma'(0), \gamma'(0))t^2+o(t^2).
\end{equation}
Denote the set of such points by $\Reg_f$, it has full measure.
When a function $f$ is $\lambda$-concave, Petrunin \cite{petrunin1996sub} proved that $\SL_f$ is a signed Radon measure. Furthermore, $\Delta^s f \ls 0$ and
\begin{equation}\label{peq:laptrace}
\Delta^{ac} f(p)=n \dashint_{\Sigma_p} H_p f(\xi,\xi) \rd \xi \ls n\cdot \lambda
\end{equation}
for almost all points $p\in M$.

We say $\SL_u \gs f \cdot \rvol$ if $\SL_{-u} \ls (-f)\cdot \rvol$. We say $\SL_u=f\cdot \rvol$ if $\SL_u \ls f\cdot \rvol$ and $\SL_u \gs f\cdot \rvol$. If $f, g\in W^{1,2}_{loc}(\Omega)$ and $\SL_g$ is a signed Radon measure, then
\begin{equation}
\begin{array}{ll}
f\SL_g(\phi)&=\int_{\Omega} \phi f \rd \SL_g\\
&=-\int_{\Omega} \langle f \nabla \phi+\phi \nabla f, \nabla g \rangle \rd \rvol
\end{array}
\end{equation}
for any $\phi \in \rLip_c(\Omega)$.

It's easy to prove the following lemma:
\begin{lem}
If $f,g, fg\in W^{1,2}_{loc}(\Omega)$ and $\SL_f ,\SL_g , \SL_{fg}$ are signed Radon measures, we have
\begin{equation}\label{eq:Lebniz}
\SL_{fg}=f\SL_g +g\SL_f +2\langle \nabla f, \nabla g \rangle \cdot \rvol.
\end{equation}
If, in addition, $f\in L^{\infty}_{loc}(\Omega)$, then we have
\begin{equation}\label{eq:chain}
\SL_{\Phi(f)}=\Phi'(f)\SL_f+\Phi''(f)|\nabla f|^2 \cdot \rvol
\end{equation}
for any $\Phi \in C^2(\bR)$.
\end{lem}

We need the following Green's formula:
\begin{lem}\label{lem:Green}
For any $p\in M$, for a.e. $R_1,R_2>0$ $(R_2>R_1)$, we have
\begin{equation}
\SL_{|p\cdot|}(B_p(R_2) \backslash \overline{B_p(R_1)}) =\cH^{n-1}(\partial B_p(R_2))-\cH^{n-1}(\partial B_p(R_1)).
\end{equation}
\end{lem}
\begin{proof}
Let $r(x)=|px|$. For $r_2>r_1>0$, denote $A_{r_1, r_2}=\{x| r_1<|px|<r_2\}$. Denote $V(t)=\rvol(B_p(t))$, by coarea formula, $V(t)=\int_0^t \cH^{n-1}(\partial B_p(s) )\rd s$. By Bishop-Gromov volume comparison theorem, $V(t)$ is locally Lipschitz. Then $V(t)$ is differentiable for a.e. $t\in (0,\infty)$ and
\begin{equation}
V'(t)=\cH^{n-1}(\partial B_p(t)).
\end{equation}
Suppose $V(t)$ is differentiable at $R_1,R_2$. Consider the cut-off functions
\[
\eta_s(x)=
\begin{cases}
\frac{|px|-R_1}{s} & \text{ if } R_1< |px| <R_1+s\\
1 &\text{ if } R_1+s\ls |px|\ls R_2-s\\
1-\frac{|px|-(R_2-s)}{s} & \text{ if } R_2-s <|px|<R_2.
\end{cases}
\]
\end{proof}
Then we have
\begin{equation}
\begin{array}{ll}
\SL_r(\eta_s)&=-\int \langle \nabla r, \nabla \eta_s \rangle\\
&=\frac1s \rvol(A_{R_2-s, R_2})-\frac1s \rvol(A_{R_1, R_1+s}).
\end{array}
\end{equation}
Let $s\to 0$, we obtain
\begin{equation}\label{ieq:limeta}
\lim_{s\to 0}\SL_r(\eta_s(x))=\cH^{n-1}(\partial B_p(R_2))-\cH^{n-1}(\partial B_p(R_1)).
\end{equation}
We claim that
\begin{equation}\label{ieq:2limeta}
\lim_{s\to 0}\SL_r(\eta_s(x))=\SL_r(A_{R_1,R_2}).
\end{equation}
Let $\SL_r=\mu^+-\mu^-$, where $\mu_{\pm}$ are Radon measures. Then
\begin{equation}
\begin{array}{ll}
|\SL_r(A_{R_1,R_2})-\SL_r(\eta_s(x))| &=|\int (\chi_{A_{R_1,R_2}}-\eta_s) \rd \mu^+ - \int (\chi_{A_{R_1,R_2}}-\eta_s) \rd \mu^-|\\
&\ls \mu^+(A_{R_1,R_1+s})+\mu^+(A_{R_2-s,R_2})\\
&+\mu^-(A_{R_1,R_1+s})+\mu^-(A_{R_2-s,R_2}).
\end{array}
\end{equation}
Let $s\to 0$, we get (\ref{ieq:2limeta}). By combining (\ref{ieq:limeta}) with (\ref{ieq:2limeta}), we finish the proof.

\textbf{Dirichlet problem}
\begin{defn}
The capacity of a set $A\subset M$ is the number
\[
C_2(A)=\inf \|u\|^2_{W^{1,2}(M)},
\]
where the infimum is taken over all $u\in W^{1,2}(M)$ such that $u\gs 1$ on $A$.
\end{defn}
It's easy to see that $C_2(\cdot)$ is countably subaddictive and $\rvol(A) \ls C_2(A)$. We say that a property regarding points in $X$ holds quasieverywhere (q.e.) if the set of points for which it fails has capacity zero.

Let $U\in M$ be a bounded domain, given a function $f\in L^2(U)$ and $g\in W^{1,2}(U)$, consider the following Dirichlet problem
\begin{equation}\label{possion}
\begin{cases}
\SL_u=f \cdot \rvol\\
u-g \in W^{1,2}_0(U).
\end{cases}
\end{equation}
If $C_2(M\backslash U)>0$, it's known that the solution exists and is unique. (See, for example, Theorem 7.12, Theorem 7.14 of \cite{cheeger1999diff} and note that if $C_2(M\backslash U)>0$, the Dirichlet Poincar\'e inequality holds, see, e.g. Corollary 5.54 of \cite{bjorn2011nonlinear}).

If $\SL_u=0$, then $u$ is called a harmonic function. If $f=0$, denote the solution of \ref{possion} by $Hg$. A Lipschitz function $g$ on $\partial U$ can be extended to a function $\tilde{g}\in \rLip(\bar{U})$ such that $g=\tilde{g}$ on $\partial U$ (see, e.g. (8.2) or (8.3) of \cite{cheeger1999diff}). By remark 7.11 and Theorem 7.14 of \cite{cheeger1999diff}, we know that $H\tilde{g}$ does not depend on the choice of extension, we define $Hg: =H \tilde{g}$. It was proved in \cite{petrunin2003harmonic} that $Hg$ is locally Lipschitz in $U$. However, it's in general not possible to have continuity up to the boundary. Nevertheless, we have the following theorem, see e.g. Theorem 10.6 of \cite{bjorn2011nonlinear}.
\begin{lem}\label{lem:continqe}
Let $U$ be a bounded domain in $M$ with $C_2(M \backslash U)>0$. Let $g:\partial U \to \bR$ be a Lipschitz function. Then for q.e. $x\in \partial U$, we have
\begin{equation}\label{continqe}
\lim_{U\ni y\to x} H g (y)=g(x).
\end{equation}
\end{lem}

We also have the following comparison principle, see e.g. Lemma 10.2 of \cite{bjorn2011nonlinear}.
\begin{lem}\label{lem:comparison}
Let $U$ be a bounded domain in $M$ with $C_2(M \backslash U)>0$. Let $g_1, g_2:\partial U \to \bR$ be two Lipschitz functions. If $g_1\ls g_2$ q.e. on $\partial U$, then $Hg_1\ls Hg_2$ in $U$.
\end{lem}

For a positive harmonic function, we have the gradient estimate, which was proved by Zhang-Zhu in \cite{zhang2012yau}, modified by Hua-Xia in \cite{hua2014anote}. Recently, Zhang-Zhu \cite{zhang2016local} have proved a sharp local Cheng-Yau gradient estimate on more general metric measure spaces.
\begin{lem}
Let $M$ be an n-dimensional Alexandrov space with curvature $\gs -K$ for some $K\gs 0$. Then there exists a constant $C=C(n)$ such that every positive harmonic function on $B_p(2R) \subset M$ satisfies
\begin{equation}\label{gradientestimate}
 |\nabla \log u| \ls C(n)(\sqrt{(n-1)K}+\frac{1}{R}) \text{ in } B_p(R).
\end{equation}
\end{lem}

\subsection{Bochner formula}
The Bochner formula for Alexandrov spaces was established in \cite{zhang2012yau}. We need the following refined version.
\begin{thm}\label{thm:sBochner}
Let $M$ be an n-dimensional Alexandrov space with curvature $\gs -K$ for some $K\gs 0$, $\partial M=\emptyset$. Suppose $u\in W^{1,2}_{loc}(M)$ and $\SL_u=f \cdot \rvol$ with $f \in W^{1,2}_{loc}(M) \cap L^{\infty}_{loc}(M)$, then $|\nabla u|^2 \in W^{1,2}_{loc}(M)$, $\Delta^s |\nabla u|^2 \gs 0$ and for $\cH^n$-a.e. $x\in M$,
\begin{equation}\label{eq:abu}
\begin{array}{ll}
[\frac12 \Delta^{ac} |\nabla u|^2 -\langle \nabla u, \nabla f \rangle +(n-1)K |\nabla u|^2-\frac{f^2}{n}] &\cdot[(1-\frac2n){\langle \nabla u, \nabla |\nabla u|^2 \rangle}^2+|\nabla u|^2 |\nabla |\nabla u|^2|^2]\\
&\gs \frac12[ |\nabla |\nabla u|^2|^2-2\frac{f}{n}\langle \nabla u ,\nabla |\nabla u|^2 \rangle]^2.
\end{array}
\end{equation}
If $\SL_u=0$, for $\cH^n$-a.e. $x\in M$,
\begin{equation}\label{bochnerweneed}
\frac12 |\nabla u|^2 \Delta^{ac} |\nabla u|^2 \gs  -(n-1)K |\nabla u|^4 +\frac{n}{4(n-1)} |\nabla |\nabla u|^2|^2.
\end{equation}
\end{thm}
To prove this theorem, we need a global Bochner formula. By \cite{petrunin2011alexmeet} and \cite{huichun2010new}, we know an n-dimensional Alexandrov space with curvature $\gs -K$ ($K\gs 0$) whose boundary is empty satisfies $RCD^*(-K(n-1),n)$ condition.
By the same trick in the proof of Theorem 3.14 of \cite{savare2013self}, we can prove the following lemma:
\begin{lem}\label{thm:gloBo}
Let $M$ be an n-dimensional Alexandrov space with curvature $\gs -K$ for some $K\gs 0$, $\partial M=\emptyset$. If $u\in W^{1,2}(M) \cap L^{\infty}(M) \cap \rLip (M)$ with $\SL_u =f\cdot \rvol$ for some $f\in W^{1,2}(M)$, then $|\nabla u|^2 \in W^{1,2}(M)$, $\Delta^s |\nabla u|^2 \gs 0$ and for $\cH^n$-a.e. $x\in M$, we have
\begin{equation}\label{gs:Bo}
\begin{array}{ll}
[\frac12 \Delta^{ac} |\nabla u|^2 -\langle \nabla u, \nabla f \rangle +(n-1)K |\nabla u|^2-\frac{f^2}{n}] &\cdot[(1-\frac2n){\langle \nabla u, \nabla |\nabla u|^2 \rangle}^2+|\nabla u|^2 |\nabla |\nabla u|^2|^2]\\
&\gs \frac12[ |\nabla |\nabla u|^2|^2-2\frac{f}{n}\langle \nabla u ,\nabla |\nabla u|^2 \rangle]^2
\end{array}
\end{equation}
\end{lem}
\begin{remark}
$\Delta^s |\nabla u|^2 \gs 0$ is by Lemma 3.2 of \cite{savare2013self}.
Note that if $n \to \infty$, then $\frac 1n \to 0$ and $\frac{n}{n-1} \to 1$, Since
\begin{equation}\label{gs:nablau}
|\nabla u|^2 |\nabla |\nabla u|^2|^2 \gs {\langle \nabla u, \nabla |\nabla u|^2 \rangle }^2,
\end{equation}
(\ref{gs:Bo}) reduces to the third inequality of Theorem 3.14 of \cite{savare2013self}.
\end{remark}
Our proof is basically along the line of the proof of Theorem 1.1 in \cite{hua2013harmonic}. We adopt their notations, denote
\[
\text{ Cutoff }=\{\psi \in \rLip_c(M): \SL \psi=\varphi \cdot \rvol \text{ for some } \varphi \in W^{1,2}(M) \cap L^{\infty}(M)\}.
\]

We need the following result on the existence of good cut-off functions. See also \cite{AMS16on,MonNarStructure}.
\begin{lem}[Propostion 2.9, \cite{hua2013harmonic}]
For any compact subset $K \subset M$, there is a $\psi \in \text{ Cutoff }$ such that $\psi=1$ in a neighborhood of $K$.
\end{lem}

\begin{lem}[Corollary 2.11, \cite{hua2013harmonic}]\label{lem:Lpsiu}
If $u\in \rLip_{loc}(M)$ with $\SL_u =f \rvol$ for some $f \in W^{1,2}_{loc}(M)\cap L^4_{loc}(M)$, then for any $\psi \in \text{ Cutoff }$, $\SL_{\psi u}=f_{\psi}\cdot \rvol$ for some $f_{\psi} \in W^{1,2}(M)$.
\end{lem}

\begin{proof}[Proof of Theorem \ref{thm:sBochner}]
Let $u\in W^{1,2}_{loc}(M)$ and $\SL_u=f \cdot \rvol$ with $f \in W^{1,2}_{loc}(M) \cap L^{\infty}_{loc}(M)$. For any Ball $B\subset M$, choose $\psi \in \mathrm{Cutoff}$ such that $\psi \equiv 1$ in a neighborhood of $\bar{B}$. Since $f\in L^{\infty}_{loc}(M)$, by \cite{Jiang2012Lip}, $u\in \rLip_{loc}(M)$. Then
\begin{equation}
\psi u \in W^{1,2}(M) \cap L^{\infty}(M) \cap \rLip(M).
\end{equation}
By Lemma \ref{lem:Lpsiu},
\begin{equation}
\SL_{\psi u} =f_{\psi} \cdot \rvol \text{ for some } f_{\psi} \in W^{1,2}(M).
\end{equation}
By Theorem \ref{thm:gloBo}, we know $|\nabla (\psi u)|^2 \in W^{1,2}(M)$ and
\begin{equation}\label{gs:abpsiu}
\begin{array}{ll}
[\frac12 \Delta^{ac} |\nabla (\psi u)|^2 -\langle \nabla (\psi u), \nabla f_{\psi} \rangle +(n-1)K |\nabla (\psi u)|^2-\frac{f_{\psi}^2}{n}] &\cdot[(1-\frac2n){\langle \nabla (\psi u), \nabla |\nabla u|^2 \rangle}^2+|\nabla u|^2 |\nabla |\nabla u|^2|^2]\\
&\gs \frac12[ |\nabla |\nabla u|^2|^2-2\frac{f_{\psi}}{n}\langle \nabla u ,\nabla |\nabla u|^2 \rangle]^2.
\end{array}
\end{equation}
Since $\psi u=u$ for any $x\in B$, we have
\begin{equation}\label{eq:psiu}
|\nabla (\psi u)|=|\nabla u|, \quad f_{\psi}=f
\end{equation}
for $x\in B$.
For any $\phi \in \rLip_c(M)$ with support in $B$, we have
\begin{equation}
\SL_{|\nabla(\psi u)|^2}(\phi)-\SL_{|\nabla u|^2}(\phi)=0.
\end{equation}
Then
\begin{equation}\label{eq:abpsiu}
\Delta^{ac}|\nabla (\psi u)|^2 =\Delta^{ac}|\nabla u|^2 \text{ for } a.e. x\in M.
\end{equation}
By combining (\ref{gs:abpsiu}), (\ref{eq:psiu}) with (\ref{eq:abpsiu}), we get (\ref{eq:abu}).
By (\ref{gs:nablau}),we have
\begin{equation}
(2-\frac2n)|\nabla u|^2 |\nabla |\nabla u|^2|^2 \gs (1-\frac2n){\langle \nabla u, \nabla |\nabla u|^2 \rangle}^2+|\nabla u|^2 |\nabla |\nabla u|^2|^2
\end{equation}
By combining this with (\ref{eq:abu}), we get (\ref{bochnerweneed}).
\end{proof}

\section{The General splitting theorem}
In this section, we prove Theorem \ref{thm:general}. We need a lemma in \cite{perelmanquasi}. We adopt some notations of this paper. Let $\Phi$ be a continuous function on $(a,b)$, $t\in (a,b)$. We write $\Phi''(t)\ls B$ if
\[
\Phi(t+\tau)\ls \Phi(t)+A\tau+\frac{B}{2}\tau^2+o(\tau^2)
\]
for some $A\in \bR$. $\Phi''(t)<\infty$ means that $\Phi''(t)\ls B$ for some $B\in \bR$. If $f$ is another continuous function on $(a,b)$, then $\Phi''\ls f$ means that $\Phi''(t)\ls f(t)$ for all $t$.
The following lemma is from 1.3 of \cite{perelmanquasi}.
\begin{lem}[1.3, \cite{perelmanquasi}]\label{lem:perconcave}
If $\Phi''(t)<\infty$ for all $t$, and $\Phi''(t)\ls f(t)+\delta$ for almost all $t$ and all $\delta>0$. Then $\Phi-F$ is concave, where $F$ is the solution of $F''=f$.
\end{lem}
We need the following Lemma:
\begin{lem}\label{2lem:mea}
Let $f:M \mapsto \bR$ be a semiconcave function with $\SL_f=c_0\cdot \rvol$ for some constant $c_0$. Let $\Phi^t$ be the $f$-gradient flow. For any Borel subset $A \subset M$, define
\begin{equation}
\Phi^{-t}(A):=\{x\in M: \Phi^t(x) \in A\}.
\end{equation}
If $\rvol(A)<\infty$, then for any $t\gs 0$,
\begin{equation}
\rvol(\Phi^{-t}(A))=\exp(-c_0t) \rvol(A).
\end{equation}
\end{lem}
This Lemma is essentially implied in the proof of 1.3. Claim of \cite{petrunin2011alexmeet}. For completeness, we present a proof here.
\begin{proof}
For any $u \in \rLip_c(M)$, $(x,t)\to u\circ \Phi^t(x)$ is locally Lipschitz. Since $M \times \bR$ is also an Alexandrov space, by Rademacher's theorem, $u\circ \Phi^t(x)$ is differentiable at $\cH^n \times \cL^1$-a.e. $y \in M \times \bR$. By Fubini theorem,
\[
I:=\{t\in [0,\infty)| u\circ \Phi^t(\cdot) \text{ is differentiable at }\cH^n-a.e. \quad x\in M\}
\]
is of full $\cL^1$-measure.
For $t\in I$,
\begin{equation}\label{geq:deriinte}
\begin{array}{ll}
\frac{d}{dt} u\circ \Phi^t(x)&=\lim_{s\to 0} \frac{(u\circ \Phi^t)(\Phi^s(x))-(u\circ\Phi^t)(x)}{s}\\
&=\langle \nabla (u\circ \Phi^t)(x),\nabla_x f \rangle
\end{array}
\end{equation}
for a.e. $x\in M$.
Since
\begin{equation}
\frac{d}{dt} u\circ \Phi^t(x)=\langle \nabla_{\Phi^t(x)}u, \nabla_{\Phi^t(x)} f \rangle
\end{equation}
and $u \in \rLip_c(M)$, $|\langle \nabla_{\Phi^t(x)}u, \nabla_{\Phi^t(x)} f \rangle| \ls C$. By (\ref{geq:deriinte}) and the dominated convergence theorem, for $t\in I$,
\[
\begin{array}{ll}
\frac{d}{dt}\int_M u\circ \Phi^t(x)\rd \rvol &=\int_M \langle \nabla (u\circ \Phi^t)(x), \nabla_x f \rangle \rd \rvol \\
&=-c_0 \int_M u\circ \Phi^t \rd \rvol.
\end{array}
\]
Denote $U(t):=\int_M u\circ \Phi^t$, then $U'(t)=-c_0 U(t)$ for a.e. $t$. Since $U(t)$ is locally Lipschitz, we have
\begin{equation}
U(t)=\exp(-c_0 t)U(0).
\end{equation}
That is,
\begin{equation}\label{2eq:mea}
\int_M u\circ \Phi^t \rd \rvol=\exp(-c_0 t) \int_M u \rd \rvol.
\end{equation}
For any ball $B=B_p(r_0)\subset M$, denote $B^r=\{x\in M: |xB| \ls r\}$. Consider the cut-off functions:
\[
u_r=
\begin{cases}
1& \text{ on } B\\
1-\frac{|xB|}{r}& \text{ on } B^r \backslash B\\
0& \text{ on } M \backslash B^r.
\end{cases}
\]

By (\ref{2eq:mea}), we have
\begin{equation}\label{geq:ur}
\int_M u_r \circ \Phi^t \rd \rvol=\exp(-c_0 t) \int_M u_r \rd \rvol.
\end{equation}

Since $u_r \to \chi_{B}$ for a.e. $x\in M$, we have
\begin{equation}\label{glim:ur}
\int_M u_r \rd \rvol \to \rvol(B).
\end{equation}
By combining (\ref{geq:ur}) with (\ref{glim:ur}), we have
\begin{equation}\label{gls:B}
\begin{array}{ll}
\rvol(\Phi^{-t}(B))&\ls \liminf_{r\to 0} \int_M u_r \circ \Phi^t \rd \rvol\\
&=\liminf_{r\to 0} \exp(-c_0t)\int_M u_r \rd \rvol\\
&=\exp(-c_0 t)\rvol(B).
\end{array}
\end{equation}

On the other hand, for $r\ls r_0$, we can choose cut-off functions $v_r$ with respect to $B=B_p(r_0)$:
\[
v_r=
\begin{cases}
1& \text{ on } B_p(r_0-r)\\
1-\frac{|xB_p(r_0-r)|}{r}& \text{ on } B \backslash B_p(r_0-r)\\
0& \text{ on } M \backslash B.
\end{cases}
\]
Then we have
\begin{equation}\label{ggs:B}
\begin{array}{ll}
\rvol(\Phi^{-t}(B))&\gs \limsup_{r\to 0} \int_M v_r \circ \Phi^t \rd \rvol\\
&=\limsup_{r\to 0} \exp(-c_0t)\int_M v_r \rd \rvol\\
&=\exp(-c_0 t)\rvol(B).
\end{array}
\end{equation}
By combining (\ref{gls:B}) with (\ref{ggs:B}), we have
\begin{equation}\label{geq:B}
\rvol(\Phi^{-t}(B))=\exp(-c_0 t)\rvol(B).
\end{equation}

Let $A$ be a Borel subset with finite volume. Denote $w_n$ the volume of unit ball $B_1(O)\subset \bR^n$. For any $\eps>0$, there exists a finite union of balls $\{B_{r_i}(p_i)\}_{i=1}^N$ such that
\begin{equation}
A \subset \cup_{i=1}^N B_{r_i}(p_i)
\end{equation}
and
\begin{equation}\label{2gs:cover}
\begin{array}{ll}
\qquad \rvol(A)&> \sum_{i=1}^N w_n r_i^n-\eps\\
&=\sum_{i=1}^N \rvol(B_{r_i}(O))-\eps\\
&>(1-\eps)\sum_{i=1}^N \rvol((B_{r_i}(p_i))-\eps.
\end{array}
\end{equation}
By combining (\ref{geq:B}) with (\ref{2gs:cover}), we have
\begin{equation}
\begin{array}{ll}
\rvol(A)&>(1-\eps)\sum_{i=1}^N \exp(c_0t) \rvol(\Phi^{-t}(B_i))-\eps\\
&\gs (1-\eps)\exp(c_0t) \rvol(\Phi^{-t}(\cup_{i=1}^N B_i))-\eps\\
&\gs (1-\eps)\exp(c_0t) \rvol(\Phi^{-t}(A))-\eps.
\end{array}
\end{equation}
By the arbitrariness of $\eps$, we have
\begin{equation}\label{ggs:A}
\rvol(A)\gs \exp(c_0 t) \rvol(\Phi^{-t}(A)).
\end{equation}

Since the Bishop-Gromov volume comparison theorem holds on Alexandrov spaces, the Vitali covering theorem follows, see for example, Theorem 1.6 of \cite{H2001lecture}. For any open subset $U \subset M$, there exist countably many disjoint balls $B_i \subset U$ such that $\rvol(U \backslash \cup_{i=1}^{\infty}B_i)=0$. By (\ref{geq:B}), we have
\begin{equation}\label{gls:U}
\begin{array}{ll}
\rvol(U)&=\rvol(\cup_{i=1}^{\infty} B_i)\\
&=\sum_{i=1}^{\infty}\rvol(B_i)\\
&=\sum_{i=1}^{\infty}\exp(c_0t) \sum_{i=1}^{\infty}\rvol (\Phi^{-t}(B_i))\\
&=\exp(c_0t)\rvol (\Phi^{-t}(\cup_{i=1}^{\infty}B_i))\\
&\ls \exp(c_0t)\rvol(\Phi^{-t}(U)).
\end{array}
\end{equation}
By combining (\ref{ggs:A}) with (\ref{gls:U}), we have
\begin{equation}\label{geq:U}
\rvol(U)=\exp(c_0t) \rvol(\Phi^{-t}(U)).
\end{equation}

Let $A\subset M$ be a Borel subset with finite volume. Then for any $\eps>0$, there exists an open subset $U \supseteq A$ such that
\begin{equation}
\rvol(U)<\rvol(A)+\eps.
\end{equation}
By (\ref{ggs:A}), we have
\begin{equation}\label{gls:UA}
\begin{array}{ll}
\rvol(\Phi^{-t}(U))-\rvol(\Phi^{-t}(A))&=\rvol(\Phi^{-t}(U\backslash A))\\
&\ls \exp(-c_0 t)\rvol(U\backslash A) \\
&< \exp(-c_0 t)\eps.
\end{array}
\end{equation}

By combining (\ref{geq:U}) with (\ref{gls:UA}), we have
\begin{equation}
\begin{array}{ll}
\rvol(A)&\ls \rvol(U)\\
&=\exp(c_0t)\rvol(\Phi^{-t}(U))\\
&< \exp(c_0t)\rvol(\Phi^{-t}(A))+\eps.
\end{array}
\end{equation}
By the arbitrariness of $\eps$, we have
\begin{equation}\label{gls:A}
\rvol(A) \ls \exp(c_0t)\rvol(\Phi^{-t}(A)).
\end{equation}
By combining (\ref{ggs:A}) with (\ref{gls:A}), we have
\[
\rvol(\Phi^{-t}(A))= \exp(-c_0t)\rvol(A),
\]
thus we complete the proof.
\end{proof}

A curve $\sigma:[0,\infty)\to M$ is called a ray if $|\sigma(s)\sigma(t)|=s-t$ for any $0\ls t<s<\infty$.

\begin{proof}[Proof of Theorem \ref{thm:general}]
Since $|\nabla_x b|=1$ for a.e. $x\in M$, $|\nabla_x f|^2-f^2(x)=0$. We will prove that $f$ is $-1$-affine, then by Theorem \ref{thm:conesplit}, $M$ splits as $M=\bR \times_{e^t} N$.

We divide our proof into four steps.

\textit{Step 1, prove that the $b$-gradient curve issuing from any point is a ray.}

Fix $R>0$, $t_0>0$. Let $\Phi^t$ be the $b$-gradient flow, consider the gradient curves $\sigma_{+,x}:[0, \infty)\to M$ of $b$ issuing from $x \in B_p(R)$. Then
\begin{equation}\label{2eq:fubini}
\int_{B_p(R)} \rd \rvol \int_0^{t_0} [1-(\nabla_{\sigma_{+,x}(t)}b)^2]\rd t =\int_0^{t_0} \int_{B_p(R)} [1-(\nabla_{\sigma_{+,x}(t)}b)^2]\rd \rvol
\end{equation}
Denote
\[
A:=\{x\in \Phi^t(B_R(p))| |\nabla_x b|\neq 1\},
\]
then $\rvol(A)=0$. Since $\SL_b=(n-1) \cdot \rvol$, by Lemma \ref{2lem:mea}, we have
\begin{equation}\label{2eq:stepmea}
\rvol(\Phi^{-t}(A))=\exp(-(n-1)t)\rvol(A)=0.
\end{equation}
By the definition of $A$,
\begin{equation}\label{2supset}
\Phi^{-t}(A)\supset \{x\in B_p(R)| |\nabla_{\sigma_{+,x}(t)}b \neq 1\}.
\end{equation}
By combining (\ref{2eq:fubini}), (\ref{2eq:stepmea}) with (\ref{2supset}), we have
\begin{equation}
\int_{B_p(R)} \rd \rvol \int_0^{t_0} [1-(\nabla_{\sigma_{+,x}(t)}b)^2]\rd t=0
\end{equation}

It follows that for a.e. $x \in B_p(R)$,
\begin{equation}\label{2eq:ae0}
\int_0^{t_0}[1-(b\circ \sigma_{+,x})'(t)] \rd t=0.
\end{equation}
Since $|\nabla_x b|=1$ a.e., $b$ is 1-Lipschitz. By (\ref{2eq:ae0}), $(b\circ \sigma_{+,x})'(t)=1$ for $\cL^1$-a.e. $t\in [0,t_0)$. It follows that, for $0\ls t_1 \ls t_2 \ls t_0$,
\begin{equation}\label{2gs:t1t2}
|\sigma_{+,x} (t_1)\sigma_{+,x} (t_2)|\gs b(\sigma_{+,x} (t_2))-b(\sigma_{+,x} (t_1))=t_2-t_1.
\end{equation}
Since $b$ is 1-Lipschitz,
\begin{equation}\label{2ls:t1t2}
|\sigma_{+,x} (t_1)\sigma_{+,x} (t_2)|\ls t_2-t_1.
\end{equation}
By combining (\ref{2gs:t1t2}) with (\ref{2ls:t1t2}), we get that for a.e. $x\in M$, $\sigma_x$ is a ray, denote this set by $M'$. For any $x \in M$, choose $M'\ni x_i \to x$, then $\sigma_{+,x_i}$ converge to the gradient curve $\sigma_{+,x}$. By Proposition \ref{prop:gradcurv}, $\sigma_{+,x}$ is a ray.

\textit{Step 2, prove that $b$ is semiconvex and the gradient curves of $-b$ and $b$ form a line.}

For any geodesic $\gamma:[0,L]\to M$ and any $t_0 \in [0,L]$, let $x=\gamma(t_0)$. For $y \in \sigma_{+,x}, y \neq x$,
\begin{equation}\label{2gs:ldiff1}
\begin{array}{ll}
b \circ \gamma(t_0+t)-b(x)&=b\circ \gamma(t)-(b(y)-|xy|)\\
&\gs |yx|-|y\gamma(t)|\\
&= \langle \uparrow_x^y, \gamma^+(t_0) \rangle t +o(t).
\end{array}
\end{equation}
Since $b(y)-b(x)=|xy|$, this means that $-|y\gamma(t)|+b(y)$ supports $b\circ \gamma(t)$ at $x$. Then $b$ is semiconvex and
\begin{equation}\label{ggs:bdiff}
d_x b(\gamma^+(t_0))\gs \langle \gamma^+(t_0), \sigma_{+,x}^+(0) \rangle.
\end{equation}
Since $-b$ is semiconcave, for $x\in M$, consider the $-b$-gradient curve $\sigma_{-,x}$ issuing from $x$. Repeating the argument as above, we can prove that $\sigma_{-,x}$ is a ray.
Denote $\sigma_x$ the curve formed by $\sigma_{-,x}$ and $\sigma_{+,x}$, let $\sigma_x(0)=x$. Since the $-b$-gradient curve issuing from $\sigma_x(t)$ for $t>0$ is a ray and geodesic doesn't branch, we know that $\sigma_x$ is a geodesic.

\textit{Step 3, prove that $b\circ \gamma$ is differentiable and estimate $(b\circ\gamma)''$.}

Let $z=\sigma_x(t)$, $t<0$. Then
\[
b(x)-b(z)=|xz|,
\]
\[
b\circ \gamma(t)-b(z)\ls |z\gamma(t)|.
\]
This means that $|z\gamma(t)|+b(z)$ supports $b\circ \gamma(t)$ at $x$. It follows that
\begin{equation}\label{gls:2diffbcosh}
(b\circ \gamma)''(t_0)\ls \sin^2 \angle (\gamma^+(t_0), \sigma_x^+(0)) \frac{\cosh|xz|}{\sinh|xz|}.
\end{equation}
and
\begin{equation}\label{gls:diffb}
\begin{array}{ll}
(b\circ \gamma)^+(t_0)&\ls -\langle \gamma^+(t_0), \sigma_x^{-}(0) \rangle\\
&=\langle \gamma^+(t_0), \sigma_x^{+}(0) \rangle.
\end{array}
\end{equation}
In (\ref{gls:2diffbcosh}), let $|xz| \to \infty$, we obtain
\begin{equation}\label{gls:2diffb}
(b\circ \gamma)''(t_0)\ls \sin^2 \angle (\gamma^+(t_0), \sigma_x^+(0)).
\end{equation}

By combining (\ref{ggs:bdiff}) with (\ref{gls:diffb}), we know that
\begin{equation}\label{geq:rdiff}
(b\circ \gamma)^+(t_0)=\langle \gamma^+(t_0), \sigma_x^{+}(0) \rangle.
\end{equation}
Repeat the argument as above, we know that if $t_0\in (0,L)$,
\begin{equation}\label{geq:ldiff}
\begin{array}{ll}
(b\circ \gamma)^-(t_0)&=\langle \gamma^-(t_0), \sigma_x^{+}(0) \rangle\\
&=-(b\circ \gamma)^+(t_0).
\end{array}
\end{equation}
By combining (\ref{geq:rdiff}) with (\ref{geq:ldiff}), we know that $b\circ \gamma(t)$ is differentiable.

Since $b$ is both semiconcave and semiconvex, so is $f=e^b$. by combining (\ref{geq:ldiff}) with (\ref{gls:2diffb}), we have
\begin{equation}
(f\circ \gamma)''(t_0) \ls f\circ \gamma(t_0),
\end{equation}
Since $f$ is semiconcave, by Lemma \ref{lem:perconcave}, we know that

\begin{equation}\label{gls:fconcave}
f \text{ is } f-\text{concave}.
\end{equation}

\textit{Step 4, prove that $f=e^b$ is $f$-affine.} Compare the proof of Lemma 4.2 of \cite{huichun2010new}.

Define the lower Hessian of $b$, $Hess_x b: T_x \to \bR$ by
\[
\underline{Hess}_x b(v,v)=\liminf_{s\to 0} \frac{b\circ \exp_x(sv)-b(x)-d_x b(v)\cdot s}{s^2/2}
\]

Since $b$ is  both semiconcave and semiconvex, it's well defined.

Recall that $Reg_b$ is the set of points $x\in M$ such that there exists Perelman's Hessian of $b$ at $x$.

Since $Reg_b$ has full measure and $\SL_b=(n-1) \cdot \rvol$, by (\ref{peq:perhess}) and (\ref{peq:laptrace}), for a.e. $x$,
\begin{equation}\label{2eq:uplow}
n\int_{\xi \in \Sigma_x}\underline{Hess}_x b(\xi,\xi)=(n-1).
\end{equation}
By (\ref{gls:2diffb}), we have
\begin{equation}\label{2gs:hess}
\underline{Hess}_x b(\xi,\xi)\ls \sin^2 \angle (\xi, \sigma_x^+(0)).
\end{equation}

By combining (\ref{2eq:uplow}) with (\ref{2gs:hess}), for a.e. $x$,
\begin{equation}\label{2eq:hess}
\underline{Hess}_x b(\xi,\xi)= \sin^2 \angle (\xi, \sigma_x^+(0)).
\end{equation}

Consider the function $\underline{u}:M \to \bR^+ \cup {0}$,
\[
\underline{u}(z)=\sup_{\xi \in \Sigma_z} |\underline{Hess}_x b(\xi,\xi)-\sin^2 \angle (\xi, \sigma_x^+(0))|.
\]

By (\ref{2eq:hess}),
\begin{equation}\label{2eq:u0}
\underline{u}=0 \text{ for a.e. } x.
\end{equation}

For any geodesic $\gamma$, by (\ref{2eq:u0}) and Segment inequality (see \cite{CC1996quanti}), there exist geodesics $\gamma_i:[0,L_i]\mapsto M$ converging to $\gamma$ uniformly such that $\int_{\gamma_i} \underline{u}(z)=0$. Then for a.e. $t$,
\begin{equation}\label{2eq:ugeo0}
\underline{u}\circ \gamma_i(t)=0.
\end{equation}
By combining (\ref{geq:ldiff}) with (\ref{2eq:ugeo0}),
for a.e. $t\in (0,L_i)$, if we denote $\xi_i^{\pm}=\gamma_i^{\pm}(t)$, then we have
\[
b \circ \gamma_i(t+s)-b\circ(t)\gs d_x b(\xi_i^+)s+\sin^2\angle (\xi_i^+, \sigma_x^+(0))+o(s^2);
\]
\[
b \circ \gamma_i(t-s)-b\circ(t)\gs -d_x b(\xi_i^-)s+\sin^2\angle (\xi_i^+, \sigma_x^+(0))+o(s^2).
\]
and $d_x b(\xi^+)=-d_x b(\xi^-)$. That is,
\begin{equation}\label{ggs:b2diff}
[b\circ \gamma_i(t)]''\gs \sin^2 \angle ((\gamma_i^+(t), (\sigma_{\gamma_i(t)})^+(0)) \text{ for a.e. } t\in (0,L_i).
\end{equation}
By combining (\ref{geq:ldiff}) with (\ref{ggs:b2diff}), we have
\begin{equation}\label{ggs:fconvex}
(f\circ \gamma_i)''(t) \gs f\circ \gamma_i(t) \text{ for a.e. } t\in (0,L_i).
\end{equation}
Since $f$ is semiconvex, by Lemma \ref{lem:perconcave}, we obtain that $f\circ \gamma_i(t)$ is $f\circ \gamma_i(t)$-convex. Since $\gamma_i$ converge to $\gamma$ uniformly, we know that
\begin{equation}\label{fconvex}
f \text{ is } f-\text{convex}.
\end{equation}

By combining (\ref{gls:fconcave}) with (\ref{ggs:fconvex}), we know that $f$ is $-1$-affine.
Since $|\nabla_x f|^2-f^2(x)=0$, by Theorem \ref{thm:conesplit}, $M$ splits as $M=\bR \times_{e^t} N$, where $N$ is an $n-1$ dimensional Alexandrov space with non-negative curvature.
\end{proof}

\section{Splitting theorem with respect to bottom of spectrum}
In this section, we will always assume that $M$ is a non-compact, n dimensional Alexandrov space with curvature $\gs -1$, $\partial M=\emptyset$. For an open subset $U\subset M$, denote $\rLip_c(U)$ the set of Lipschitz functions with compact support in $U$. The bottom of the $L^2$ spectrum of the Laplacian on $M$ can be characterized by
\[
\lambda_0(M)=\inf_{f\in \rLip_c(M)} \frac{\int_M |\nabla f|^2 \rd \rvol}{\int_M f^2 \rd \rvol}.
\]

Now fix a ball $B_p(R_0)$, from now on, we say $E$ is an end of $M$, we mean $E$ is an unbounded connected component of $M\backslash B_p(R_0)$.
Let $E$ be an end of $M$. The bottom of the $L^2$ spectrum of the Laplacian on $E$ satisfying Dirichlet boundary condition on $\partial E$ can be characterized by
\[
\lambda_0(E)=\inf_{f\in \rLip_c(E)} \frac{\int_E |\nabla f|^2 \rd \rvol}{\int_E f^2 \rd \rvol}.
\]
It's easy to see that $\lambda_0(E) \gs \lambda_0(M)$.

We adopt some notations of \cite{LiWang2001complete1}. If $E$ is an end of $M$, denote $E(R)=E \cap B_p(R)$ and $\partial E(R)=E \cap \partial B_p(R)$. Denote $V_E(\infty)=\rvol(E)$, $V_E(R)=\rvol (E\cap B_p(R))$.
Now let $R_0<R_1<R_2<...\to \infty$,
Consider the harmonic functions:
\[
\SL_{h_R}=0 \text{ on } E(R_i),
\]
\[
h_R=1 \text{ on } \partial E,
\]
and
\[
h_R=0 \text{ on } \partial E(R_i).
\]
By the maximum principle, $0\ls h_{R_i}\ls 1$. By gradient estimate (\ref{gradientestimate}), on any compact subset of $E$, $h_{R_i}$ is equi-continuous for sufficiently large $R_i$. By Arzela-Ascoli's theorem, there exists a subsequence converging locally uniformly to a Lipschitz function $h$ defined on $E$, $0\ls h\ls 1$. By Lemma \ref{plem:conver}, $\SL_h=0$, $h$ is harmonic. Note that $h$ may be a constant.

\begin{lem}\label{lem:increasing}
If $R_i\ls R_j$, then $h_{R_i}\ls h_{R_j}$ on $E(R_i)$.
\end{lem}
\begin{proof}
By Lemma \ref{lem:continqe}, for any $k$, for q.e. $x\in \partial E$,
\begin{equation}\label{eq:qe1}
\lim_{E\ni y\to x} h_{R_k}(y)=1.
\end{equation}
Then for q.e. $x\in \partial E$,
\begin{equation}
\lim_{E\ni y\to x} (h_{R_j}-h_{R_i})(x)=0.
\end{equation}
Note that
\begin{equation}
h_{R_j}|_{\partial E(R_i)}\gs 0=h_{R_i}.
\end{equation}
By Lemma \ref{lem:comparison}, $h_{R_i}\gs h_{R_j}$ on $E(R_i)$.
\end{proof}

\begin{defn}\label{defn:Alexnonpara}
An end $E$ is said non-parabolic if the sequence of harmonic functions $h_{R_i}$ subconverge to a non-constant harmonic function $h$. Otherwise, it's said parabolic.
\end{defn}

\begin{lem}
If $E$ is a non-parabolic end, then $\inf_E h=0$.
\end{lem}
\begin{proof}
Let $c=\inf_E h$, then $0\ls c<1$. Consider $\tih:=\frac{h-c}{1-c}$, then
\begin{equation}\label{ls:thh}
0\ls \tih \ls h \ls 1.
\end{equation}
Since $h_{R_k}\ls h \ls 1$ for any $k$, by (\ref{eq:qe1}), we have
\begin{equation}
\lim_{E\ni y\to x} h(y)=1 \text{ for } q.e. x \in \partial E.
\end{equation}
It follows that
\begin{equation}\label{eq:thqe1}
\lim_{E\ni y\to x} \tih(y)=1 \text{ for } q.e. x \in \partial E.
\end{equation}
Then for any $k$,
\begin{equation}
\lim_{E\ni y\to x} (\tih-h_{R_k})(y)=0 \text{ for } q.e. x \in \partial E.
\end{equation}
Since
\begin{equation}
\tih|_{\partial E(R_k)}\gs 0=h_{R_k}.
\end{equation}
By Lemma \ref{lem:comparison}, $\tih\ls h_{R_k}$ on $E(R_k)$. It follows that
\begin{equation}\label{gs:thh}
\tih \gs h.
\end{equation}
By combining (\ref{ls:thh}) with (\ref{gs:thh}), we have $\tih=h$, then $c=0$.
\end{proof}

\begin{remark}
Let $E$ be a non-parabolic end.
If $h_{R_i}$ subconverge to another harmonic function $h'$, since $h'\gs h_{R_k}$ for any $k$, we have $h'\gs h$. Similarly, we can get $h\gs h'$, then $h'=h$.
If $\rho_i \to \infty$, suppose $h_{\rho_i}$ subconverge to a harmonic function $h''$ defined on $E$. Repeat the above argument, we can get $h''=h$.
If $E$ is a parabolic end, then $h\equiv 1$.
\end{remark}

\begin{lem}\label{lem:bddhar}
Suppose $M$ has at least two non-parabolic ends, then there exists a non-constant, bounded harmonic function defined on $M$.
\end{lem}
The proof of this proposition is similar to the case of Riemannian manifolds, see \cite{LiTam1992harmonic}. We include a proof here.
\begin{proof}
Suppose $R_0>0$ is sufficiently large so that $M \backslash B_p(R_0)$ has at least two disjoint non-parabolic ends $E_1$ and $E_2$. Choose an increasing sequence $R_i \to \infty$ such that $R_1>R_0$, let $f_{R_i}$ be the solution of
\[
\SL_{f_{R_i}}=0 \text{ on } B_p(R_i),
\]
\[
f_{R_i}=1 \text{ on } \partial E_1(R_i),
\]
and
\[
f_{R_i}=0 \text{ on } \partial B_p(R_i) \backslash E_1.
\]
Clearly, $\partial E_2(R_i) \subset \partial B_p(R_i) \backslash E_1$. Then $f_{R_i}$ subconverge to a harmonic function $f$ satisfying $0\ls f\ls 1$. Next, we prove that $f$ is not a constant. For $k=1,2$, let $h_{k,R_i}$ be the harmonic functions on $E_k(R_i)$ such that
\[
h_{k, R_i}|_{\partial E_k}=1, \quad h_{k, R_i}|_{\partial E_k(R)}=0.
\]
Suppose $h_{k,R_i}$ subconverge to a harmonic function $h_k$ defined on $E_k$. By lemma \ref{lem:continqe}, for q.e. $x\in \partial E_2$,
\begin{equation}
\lim_{E_2\ni y \to x} h_{2,R_i}=1 \gs f_{R_i}(x).
\end{equation}
Note that
\begin{equation}
h_{2,R_i}=0=f_{R_i} \text{ on } \partial E_2(R_i),
\end{equation}
By Lemma \ref{lem:comparison}, we have
\[
h_{2,R_i} \gs f_{R_i} \text{ on } E_2(R_i).
\]
It follows that $h_2 \gs f$ on $E_2$. Since $\inf_{E_2} h_2=0$,
\begin{equation}\label{eq:E20}
\inf_{E_2} f=0.
\end{equation}
By repeating the above argument, we can prove that
\[
1-h_{1,R_i} \ls f_{R_i} \text{ on } E_2(R_i).
\]
Then
\[
1-h_1 \ls f \text{ on }E_1.
\]
Since $\inf_{E_1} h_1 =0$, we know that
\begin{equation}\label{eq:E11}
\sup_{E_1} f=1.
\end{equation}
By combining (\ref{eq:E20}) with (\ref{eq:E11}), we know that $f$ is non-constant.
\end{proof}

Following the argument in the proof of Theorem 22.1 of \cite{Li2012Geometric}, we can get the following decay estimate.
\begin{lem}\label{lem:decayest}
Let $M$ be an n-dimensional Alexandrov space with curvature $\gs -K$ for some $K\gs 0$, $\partial M=\emptyset$. Suppose $E$ is an end of $M$ with respect to $B_p(R_0)$ such that $\lambda_0(E) >0$. Let $f$ be a non-negative function defined on $E$ satisfying $\SL_f \gs 0$. If $f$ satisfies the growth condition
\begin{equation}\label{eq:growth}
\int_{E(R)} f^2 \exp (-2\sqrt{\lambda_0(E)}r)=o(R)
\end{equation}
as $R \to \infty$, then it must satisfies the decay estimate
\begin{equation}
\int_{E(R+1)\backslash E(R)} f^2 \ls C(1+(R-R_0)^{-1})\exp(-2\sqrt{\lambda_0(E)}R)\int_{E(R_0+1)\backslash E(R_0)} \exp(2\sqrt{\lambda_0(E)}r) f^2.
\end{equation}
for some constant $C>0$ depending on $\lambda_0(E)$ and for all $R \gs 2(R_0+1)$.
\end{lem}

Suppose $E_1$ is an end of $M$. Let $R_i \to \infty$ be an increasing sequence, consider the harmonic functions
\[
\SL_{f_{R_i}}=0 \text{ on } B_p(R_i),
\]
\[
f_{R_i}=1 \text{ on } \partial E_1(R_i),
\]
and
\[
f_{R_i}=0 \text{ on } \partial B_p(R_i) \backslash E_1.
\]
Then $f_{R_i}$ subconverge to a harmonic function $f$ defined on $M$. Note that $f$ may be a constant.
We can get the the following decay estimate for $f$. See Corollary 22.3 of \cite{Li2012Geometric} and Lemma 1.1 of \cite{LiWang2001complete1}.
\begin{lem}\label{lem:f-1}
Suppose $E_1$ is an end of $M$, $f$ is the harmonic function constructed above. If $\lambda_0(E_1)>0$, then
\begin{equation}\label{ls:f-1}
\int_{E_1(R+1)\backslash E_1(R)}(f-1)^2 \ls C_1 \exp (-2R \sqrt{\lambda_0(E_1)})
\end{equation}
for some constant $C_1>0$ depending on $f$, $\lambda_0(E_1)$ and $n$.
If $E$ is another end with $\lambda_0(E)>0$, then
\begin{equation}\label{ls:ef-1}
\int_{E(R+1)\backslash E(R)} f^2 \ls C \exp (-2R \sqrt{\lambda_0(E)})
\end{equation}
for some constant $C>0$ depending on $f$, $\lambda_0(E)$ and $n$.
\end{lem}
\begin{proof}
Consider the functions
\[
\tf_{R_i}=
\begin{cases}
f_{R_i} &\text{ on }E_1(R_i)\\
1 & \text{ on } E_1 \backslash B_p(R_i).
\end{cases}
\]
Let $g_{R_i}$ be a Lipschitz function defined on $\overline{B_p(R_i)}$ such that
\[
g_{R_i}|_{\partial E_1(R_i)}=1, \quad g_{R_i}|_{\partial B_p(R_i) \backslash E_1}=0.
\]
Let
\[
\tg_{R_i}=
\begin{cases}
g_{R_i} &\text{ on }E_1(R_i)\\
1 & \text{ on } E_1 \backslash B_p(R_i).
\end{cases}
\]
then $\tg_{R_i}$ is Lipschitz. Since $f_{R_i}-g_{R_i} \in W^{1,2}_0(B_p(R))$,
$1-\tf_{R_i} \in W^{1,2}(E_1)$. It's easy to see that $1-\tf_{R_i}$ satisfies the growth condition (\ref{eq:growth}). By Lemma 7.13 of \cite{bjorn2011nonlinear}, we have
\begin{equation}
\SL_{1-\tf_{R_i}} \gs 0.
\end{equation}
By Lemma \ref{lem:decayest}, we can get
\begin{equation}\label{ls:tfRi}
\int_{E_1(R+1)\backslash E_1(R)}(\tf_{R_i}-1)^2 \ls C_1 \exp (-2R \sqrt{\lambda_0(E_1)})
\end{equation}
for some constant $C_1>0$ depending on $f$, $\lambda_0(E_1)$ and $n$. Note that $1-\tf_{R_i}=1-f_{R_i}$ on $E_1(R_i)$ and vanishes on $E_1\backslash B_p(R_i)$. Then if we replace $\tf_{R_i}-1$ by $f_{R_i}-1$, (\ref{ls:tfRi}) still holds. By letting $R_i \to \infty$, we get (\ref{ls:f-1}). Similarly, we can get (\ref{ls:ef-1}).
\end{proof}

Following the argument in the proof of Lemma 1.2 of \cite{LiWang2001complete1}, we can get
\begin{lem}\label{lem:integral}
If $E$ is an end of $M$ with $\lambda_0(E)>0$, the harmonic function $f$ in Lemma \ref{lem:f-1} satisfies
\[
\int_{E(R)} \exp(2\sqrt{\lambda_0(E)}r)|\nabla f|^2 \ls CR
\]
for $R$ sufficiently large.
\end{lem}

Li-Wang \cite{LiWang2001complete1} proved sharp volume growth/decay rates for an end $E$ with $\lambda_0(E)>0$, see Theorem 1.4 of \cite{LiWang2001complete1}. This has been generalized by Buckley-Koskela \cite{BK2006ends} to proper pointed metric measure spaces, which include Alexandrov spaces. To state the estimate, denote by $V_E(R)$ the volume of the set $E(R)$. The volume of the end $E$ will be denoted by $V_{E}(\infty)$.
\begin{lem}\label{lem:volume}
Let $E$ be an end of $M$ with $\lambda_0(E)>0$.

(1) If $E$ is a parabolic end, then $E$ must have exponential volume decay given by
\begin{equation}\label{ls:paravolume}
V_E(\infty)-V_E(R) \ls C \exp(-2 \sqrt{\lambda_0(E)}R)
\end{equation}
for some constant $C>0$ depending on the end $E$.

(2) If $E$ is a non-parabolic end, then $E$ must have exponential volume growth given by
\begin{equation}\label{gs:nonparavolume}
V_E(R)\gs C\exp(2\sqrt{\lambda_0(E)}R)
\end{equation}
for some constant depending on the end $E$.
\end{lem}
\begin{remark}
For a parabolic end $E$ with $\lambda_0(E)>0$, following the argument in the proof of (1) of Theorem 1.4 in \cite{LiWang2001complete1}, we can prove the estimate (\ref{ls:paravolume}). Buckley-Koskela proved that if an end $E$ satisfies $\lambda_0(E)>0$, then it either has volume decay as (\ref{ls:paravolume}) or has volume growth as (\ref{gs:nonparavolume}). So non-parabolic ends must satisfies (\ref{gs:nonparavolume}).
\end{remark}

The following theorem, when restricted to Riemannian manifolds, is a particular case of Theorem 2.1 of \cite{LiWang2001complete1}.
\begin{lem}\label{lem:onenon}
Suppose $n\gs 3$. Let $M$ be a non-compact, n dimensional Alexandrov space with curvature $\gs -1$, $\partial M=\emptyset$. If $\lambda_0>n-2$, then $M$ has only one end with infinite volume.
\end{lem}
\begin{proof}
We argue by contradiction. Suppose $M$ has two ends $E_1,E_2$ with infinite volume. By Lemma \ref{lem:volume}, we know they are non-parabolic. By Lemma \ref{lem:bddhar}, there exists a non-constant, bounded harmonic function $f$ defined on $M$. Let $\psi=|\nabla f|^2$, by Theorem \ref{thm:sBochner}, $\psi \in W^{1,2}_{loc}(M) \cap L^{\infty}_{loc}(M)$, $\Delta^s \psi \gs 0$ and
\begin{equation}\label{gs:psi}
\psi \Delta^{ac} \psi +2(n-1) \psi^2-\frac{n}{2(n-1)}|\nabla \psi|^2\gs 0 \text{ for } \ a.e. x.
\end{equation}
Denote $g=|\nabla f|^{\frac{n-2}{n-1}}$, by the following Lemma \ref{lem:Lg}, $g\in W^{1,2}_{loc}(M)$ and
\[
\SL_g \gs -(n-2) g\cdot \rvol.
\]
By Lemma \ref{lem:integral} and following the argument from line 11 on page 520 to line 7 on page 521 of \cite{LiWang2001complete1}, we can prove that
\[
\int_{B_p(2R) \backslash B_p(R)} g^2 \ls CR.
\]
Following the argument from line 8 on page 521 to line 9 on page 522 of \cite{LiWang2001complete1}, we can find non-negative functions $\phi_R \in \rLip_c(M)$ such that
\[
\int_M |\nabla (\phi_R f)|^2 \ls (n-2) \int_M \phi_R^2 g^2 +\int_M |\nabla \phi_R|^2 g^2
\]
and
\[
\int_M |\nabla \phi_R|^2 g^2 \ls CR^{-2} \int_{B_p(2R) \backslash B_p(R)} g^2 \to 0.
\]
It follows that $\lambda_0 \ls n-2$, contradiction! Hence we complete the proof.
\end{proof}

\begin{lem}\label{lem:Lg}
Let $g=|\nabla f|^{\frac{n-2}{n-1}}$, then $g\in W^{1,2}_{loc}(M)$ and
\begin{equation}\label{Lggs}
\SL_g \gs -(n-2) g\cdot \rvol.
\end{equation}
\end{lem}
\begin{proof}
Let $\psi=|\nabla f|^2$, by Theorem \ref{thm:sBochner}, $\psi \in W^{1,2}_{loc}(M) \cap L^{\infty}_{loc}(M)$, $\Delta^s \psi \gs 0$ and
\begin{equation}\label{gs:psi}
\psi \Delta^{ac} \psi +2(n-1) \psi^2-\frac{n}{2(n-1)}|\nabla \psi|^2\gs 0.
\end{equation}
Following the argument in the proof of Lemma 4.12 of \cite{PRG2008}, we can prove that for $p>\frac{1-\frac{n}{2(n-1)}}{2}=\frac{n-2}{4(n-1)}$, $\psi^p\in W^{1,2}_{loc}(M)$ and furthermore,
\begin{equation}\label{psiepsto}
\SL_{(\psi^2+\eps)^{\frac{p-1}{2}}\psi}(\varphi) \to \SL_{\psi^p}(\varphi).
\end{equation}
for $\varphi \in \rLip_c(M)$. By (\ref{eq:chain}), we have
\begin{equation}\label{eq:Lpsi}
\SL_{(\psi^2+\eps)^{\frac{p-1}{2}}\psi}
=(\psi^2+\eps)^{\frac{p-3}{2}}(p\psi^2+\eps)\SL_{\psi}+[(\psi^2+\eps)^{\frac{p-5}{2}}[p(p-1)\psi^3+3\eps(p-1)\psi^3]|\nabla \psi|^2]\cdot \rvol.
\end{equation}
Now let $p=\frac{(n-2)}{2(n-1)}$, then
\begin{equation}\label{sinpsip}
\Delta^s \psi^p \gs 0.
\end{equation}
By (\ref{gs:psi}) and (\ref{eq:Lpsi}), for a.e. $x\in M$,
\begin{equation}\label{Lpsieps}
\begin{array}{ll}
\SL_((\psi^2+\eps)^{\frac{p-1}{2}}\psi)
&\gs \Delta^{ac}((\psi^2+\eps)^{\frac{p-1}{2}}\psi)\\
&\gs [\eps(\psi^2+\eps)^{\frac{p-3}{2}}-\frac{4(n-1)^2}{n}(\psi^2+\eps)^{\frac{p-5}{2}}[p(p-1)\psi^5+3\eps(p-1)\psi^5]]\cdot \rvol\\
&\gs [-\frac{4(n-1)^2}{n}(\psi^2+\eps)^{\frac{p-5}{2}}[p(p-1)\psi^5+3\eps(p-1)\psi^5]]\cdot \rvol.
\end{array}
\end{equation}
Let $\eps \to 0$, by (\ref{psiepsto}), we have
\begin{equation}
\SL_{\psi^p} \gs -(n-2)\psi^p \cdot \rvol
\end{equation}
Note that $\psi^p=|\nabla f|^{\frac{n-2}{n-1}}=g$, thus we get (\ref{Lggs}).
\begin{remark}
In Lemma 4.12 of \cite{PRG2008}, it's assumed that $\psi \in \rLip_{loc}(M)$. We find that for $\psi \in W^{1,2}_{loc}(M)\cap L^{\infty}_{loc}(M)$, the lemma still holds.
\end{remark}
\end{proof}
Next, we prove Theorem \ref{thm:eigensplitting}. We follow the argument in the proof of Theorem 1.1 of \cite{JLW2009ends}.

\begin{proof}[Proof of Theorem \ref{thm:eigensplitting}]
If (1) doesn't hold, then $M$ has at least two ends. For $n\gs 4$, $\frac{(n-1)^2}{4}>n-2$. By Lemma \ref{lem:onenon}, $M$ has at most one non-parabolic end, then $M$ has at least one parabolic end $E$. Let $\gamma: [0,\infty) \to M$ be a ray with $\gamma(0)=p$, $\gamma(t) \to E(\infty)$, where $E(\infty)$ denotes the infinity of the end $E$. Consider the Busemann function (note that it's different from the common form) $b:M \to \bR$,
\[
b(x)=\lim_{t \to \infty} |x\gamma(t)|-t.
\]
For any $x\in M$, choose a geodesic $\gamma_{x,i}$ connecting $x$ with $\gamma(i)$, let $\gamma_{x,i}(0)=x$. Then there exists a subsequence $\gamma_{x,i_j}$ of $\gamma_{x,i}$ converging to a ray $\gamma_x$. Note that it may not be unique. For any $t_0>0$, we have
\begin{equation}\label{eq:buseray}
\begin{array}{ll}
b\circ \gamma_x(t_0)-b(x)&=\lim_{j \to \infty}|\gamma(i_j)\gamma_x(t_0)|-|\gamma(i_j)x|\\
&=\lim_{j \to \infty}|\gamma(i_j)\gamma_x(t_0)|-|\gamma(i_j)\gamma_{x,i_j}(t_0)|-|x\gamma_{x,i_j}(t_0)|\\
&=\lim_{j \to \infty}(|\gamma(i_j)\gamma_x(t_0)|-|\gamma(i_j)\gamma_{x,i_j}(t_0)|)-t_0\\
&=-t_0.
\end{array}
\end{equation}
The last inequality holds since
\[
||\gamma(i_j)\gamma_x(t_0)|-|\gamma(i_j)\gamma_{x,i_j}(t_0)||\ls |\gamma_x(t_0)\gamma_{x,i_j}(t_0)|\to 0.
\]
It follows that for any $x \in M$,
\begin{equation}\label{eq:Lipb}
\rLip b(x)=1.
\end{equation}
By a similar argument in Step 3 in the proof of Theorem \ref{thm:general}, we can prove that $b$ is semiconcave and for any geodesic $\sigma$ (let $x=\sigma(t)$),
\begin{equation}
(b\circ \sigma)''(t) \ls \sin^2 \angle (\gamma_x^+(0), \sigma'(t)).
\end{equation}
By (\ref{peq:laptrace}), for a.e. $x\in M$,
\begin{equation}
\Delta^{ac} b(x) \ls n-1.
\end{equation}
Since $\Delta^s b \ls 0$, we have
\begin{equation}\label{els:Lbn-1}
\SL_b \ls (n-1)\cdot \rvol.
\end{equation}
By (\ref{eq:Lipb}) and Lemma \ref{lem:graequal}, we have
\begin{equation}\label{eq:gra1}
|\nabla_x b|=1 \text{ for } a.e. x\in M.
\end{equation}

Denote
\begin{equation}\label{eq:u}
u=\exp(-\frac{n-1}{2}b),
\end{equation}
we have
\begin{equation}\label{gs:Lu}
\begin{array}{ll}
\SL_u&=-\frac{n-1}{2} \exp(-\frac{n-1}{2}b) \SL_b+[\frac{(n-1)^2}{4}  \exp(-\frac{n-1}{2}b)]\cdot \rvol\\
&\gs [-\frac{(n-1)^2}{4} u] \cdot \rvol.
\end{array}
\end{equation}
We will prove that $\SL_u =-\frac{(n-1)^2}{4} u\cdot \rvol$.
For any non-negative function $\phi\in Lip_c(M)$,
\[
\begin{array}{ll}
\int_M |\nabla(\phi u)|^2&=-\SL_u (\phi^2 u)+\int_M |\nabla \phi|^2 u^2\\
&=\frac{(n-1)^2}{4}\int_M \phi^2 u^2 +\int_M |\nabla \phi|^2 u^2-\SL_u(\phi^2 u)-\frac{(n-1)^2}{4}\int_M \phi^2 u^2.
\end{array}
\]
Since
\begin{equation}
\int_M |\nabla (\phi u)|^2 \gs \frac{(n-1)^2}{4}\int_M \phi^2 u^2,
\end{equation}
We have
\begin{equation}\label{ls:eigeninte}
\SL_u(\phi^2 u)+\frac{(n-1)^2}{4}\int_M \phi^2 u^2 \ls \int_M |\nabla \phi|^2 u^2.
\end{equation}
Following the argument from line 20 on page 5 to line 26 one page 6 of \cite{JLW2009ends}, there exist $\phi_R \in Lip_c(M)$ such that:
For $x\in E$,
\[
\phi_R(x)=
\begin{cases}
1 & \text{ if }|px| \ls R\\
\frac{2R-|px|}{R} & \text{ if } R \ls |px| \ls 2R\\
0 & \text{ if } |px|\gs 2R.
\end{cases}
\]
and
\begin{equation}\label{ls:painte}
\int_E |\nabla \phi_R|^2 u^2 \ls C R^{-1}.
\end{equation}

For $x \in M \backslash E$,
\[
\phi_R(x)=
\begin{cases}
1 & \text{ if } b(x)\ls R\\
\frac{2R-b(x)}{R} & \text{ if } R \ls b(x) \ls 2R\\
0 & \text{ if } b(x) \gs 2R.
\end{cases}
\]
and
\begin{equation}\label{ls:noninte}
\int_{M \backslash E} |\nabla \phi_R|^2 u^2 =R^{-2} \int_{\bar{B}(R,2R)} \exp(-(n-1)b),
\end{equation}
where
\[
\bar{B}(R,2R)=\{x\in M\backslash E| R \ls b(x)<2R \}.
\]
We now claim that the volume of $\bar{B}(R,R+1)$, denoted by $\bar{V}(R,R+1)$, is bounded by $C\exp((n-1)R)$ for sufficiently large $R$. By (\ref{els:Lbn-1}), we have
\[
(n-1) \bar{V}(R_1,R_2)\gs \SL_b(\bar{B}(R_1,R_2)).
\]
Following the argument of the proof of Lemma \ref{lem:Green}, we can prove that
\[
\SL_b(\bar{B}(R_1,R_2))=\bar{A}(R_2)-\bar{A}(R_1)
\]
for a.e. $R_1, R_2$, where $\bar{A}(R)$ denotes the $n-1$ dimensional Hausdorff meansure of $\{x\in M \backslash E| b(x)=R\}$. Note that $\bar{V}(R_0,t)$ is locally Lipschitz, following the same argument from line 1 to line 10 on page 7 of \cite{JLW2009ends}, we prove the claim and get
\begin{equation}\label{ls:noninte2}
\int_{\bar{B}(R,2R)} \exp(-(n-1)b) \ls CR.
\end{equation}
By combining (\ref{ls:painte}), (\ref{ls:noninte}) and (\ref{ls:noninte2}), we have
\begin{equation}\label{ls:Minte}
\int_M|\nabla \phi_R|^2 u^2 \to 0.
\end{equation}
By combining (\ref{ls:eigeninte}) with (\ref{ls:Minte}), we have
\begin{equation}\label{inteto0}
\int_M \phi^2_R u \rd \Delta^s u+ \int_M (\Delta^{ac}u +\frac{(n-1)^2}{4}u) \phi^2_R u \rd \rvol \to 0.
\end{equation}
Note that by (\ref{gs:Lu}), the measure $\Delta^s u$ is non-negative and $\Delta^{ac} u \gs -\frac{(n-1)^2}{4}u$ almost everywhere. Then the first term and the second term of the left hand side of (\ref{inteto0}) are non-negative. It follows that
\begin{equation}
\int_M (\Delta^{ac}u +\frac{(n-1)^2}{4}u) \phi^2_R u \rd \rvol \to 0,
\end{equation}
and
\begin{equation}\label{sinto0}
\int_M \phi^2_R u \rd \Delta^s u \to 0.
\end{equation}
Thus
\begin{equation}\label{eq:ac0}
\Delta^{ac}u +\frac{(n-1)^2}{4}u=0 \text{ for } \cH^n a.e. x \in M.
\end{equation}
We claim that
\begin{equation}\label{eq:sin0}
\Delta^s u(M)=0.
\end{equation}
Otherwise, there exist $R_1>0$ and $C>0$, such that
\[
\Delta^s u(\bar{B}_p(R_1))>C>0.
\]
Then
\[
\int_{B_p(R)} \phi_R^2 u \rd \Delta^s u \gs C \min_{x\in B_p(R_1)} u
\]
for $R$ sufficiently large. This contradicts to (\ref{sinto0}), thus (\ref{eq:sin0}) holds. By combining this with (\ref{eq:ac0}), we have
\begin{equation}\label{eq:Lu}
\SL_u=-\frac{(n-1)^2}{4}u \cdot \rvol.
\end{equation}
By combining this with (\ref{eq:u}), we have
\begin{equation}
\SL_b =(n-1) \cdot  \rvol.
\end{equation}
Since $b$ is semiconcave and $|\nabla_x b|=1$ for a.e. $x \in M$, by Theorem \ref{thm:general}, $M$ splits and $M=\bR \times_{e^t} N$, where $N$ is an n-1 dimensional Alexandrov space with non-negative curvature. Since $M$ has at least two ends, by the argument in the proof of Lemma 9.5 of \cite{PRG2008}, we know that $N$ is compact.
\end{proof}

\begin{proof}[Proof of Theorem \ref{thm:An23}]
If $M$ has one finite volume end (i.e. parabolic end), the proof is the same as above.
\end{proof}
\begin{remark}\label{rem:3weaker}
For $n=3$, our result is weaker than Theorem \ref{thm:ngs3}.  If $M$ has two infinite volume end, $\lambda_0(M)=1$, we don't know whether $M$ splits as case (3) of Theorem \ref{thm:ngs3}. Since in the proof of Lemma \ref{lem:onenon}, we don't know whether we can get rigidity from $\SL_g=-(n-2)g \cdot \rvol$.
\end{remark}
\begin{remark}\label{rem:RCD}
Let $(X,d,m)$ be a complete, separable metric measure space satisfying the Riemannian curvature-dimension condition $RCD^*(-(N-1),N)$. For the Dirichlet problem, Lemma \ref{lem:continqe} and \ref{lem:comparison} holds for metric measure spaces with a doubling measure and satisfying a $(1,p)$ Poincar\'e inequality for $p>1$. $RCD^*(K,N)$ spaces are included. So we can define parabolic ends and non-parabolic ends similarly. Suppose $\lambda_0(X)>0$, by Theorem 0.1 of \cite{BK2006ends}, the volume of these ends satisfies exponential growth/decay estimates. Decay estimates for harmonic functions also holds since test functions are compositions of distance functions. So we can prove an analogue of Lemma \ref{lem:onenon}.

If $X$ has a finite volume end $E$, let $b$ be the Busemann function with respect to the ray to the infinity of $E$. Following Gigli's argument in \cite{gigli2013split,gigli2013overview}, we may prove that $\SL_b \ls (n-1) \cdot \rvol$. Following the proof of Theorem \ref{thm:eigensplitting}, we may prove that $\SL_b=(n-1)\cdot m$ and the minimal relaxed gradient $|\nabla b|_w=1$ for $m$-a.e.$x\in X$. However, we don't know whether $b$ is semiconcave, since $(X,d,m)$ has only "Ricci curvature bounded below". For any $x\in X$, we don't know whether the gradient curve of $b$ exists. So it seems to me that our argument can't be generalized to $RCD^*(K,N)$ directly. However, analogue theorems may hold on $RCD^*(K,N)$ spaces.
\end{remark}

\section{Splitting theorem with respect to volume entropy}
In this section, we always suppose that $M$ is a compact, n-dimensional Alexandrov space with curvature $\gs -1$. Since Alexandrov space is locally contradictable, the universal cover $\pi: \tilde{M} \mapsto M$ exists. We are concerned with the volume entropy $h$ defined by
\begin{equation}
h(M)=\lim_{r \to \infty} \frac{\ln \rvol(B_{\tM}(x,r))}{r},
\end{equation}
By the same argument in \cite{M1979ventropy}, the limit exists and is independent of the center $x\in \tM$. By the volume comparison theorem, we know that $h\ls n-1$.

\begin{proof}[Proof of Theorem \ref{thm:volent}]
$\Psi(\delta)$ means that when $\delta \to 0$, $\Psi \to 0$.
First, we follow the approach of \cite{Liu2011short} to  construct a Busemann function $u$ on $\tM$ and show that $\SL_u=(n-1)\cdot \rvol$.
Now take a fixed $R> 50 diam M$. Pick a point $O\in \tM$ and define $r(x)=|Ox|$. Following the same argument in the proof of Claim 1 of \cite{Liu2011short}, we can prove that: there exists a sequence $r_i \to \infty$ such that
\begin{equation}\label{heq:infty}
\frac{\cH^{n-1}(\partial B(O, r_i+50 R))}{\cH^{n-1}(\partial B(O, r_i-50 R))}=\exp[100(n-1)R-\Psi(\frac1i)]
\end{equation}
Now define
\[
A_i=\{x\in \tM|r_i-50R \ls r(x) \ls r_i+50R \}.
\]
By Lemma \ref{lem:Green}, without loss of generality, we can assume that
\begin{equation}\label{heq:green}
\SL_r(A_i)=\cH^{n-1}(\partial B(O, r_i+50 R))-\cH^{n-1}(\partial B(O, r_i-50 R)).
\end{equation}
By the relative volume comparison theorem,
\begin{equation}\label{hls:BG}
\begin{array}{ll}
\rvol(A_i)&\ls \cH^{n-1}(\partial B(O, r_i-50 R))\int_{r_i-50R}^{r_i+50R} (\frac{\sinh t}{\sinh (r_i-50R)})^{n-1} \\
&\ls \int_0^{100R} \exp((n-1)t) \rd t \ \cH^{n-1}(\partial B(O, r_i-50 R))\\
&\ls \frac{\exp(100(n-1)R)-1}{n-1} \cH^{n-1}(\partial B(O, r_i-50 R)).
\end{array}
\end{equation}
By combining (\ref{heq:infty}), (\ref{heq:green}) with (\ref{hls:BG}), we have
\begin{equation}
\frac{\SL_r(A_i)}{\rvol(A_i)} \gs n-1-\Psi(\frac1i).
\end{equation}

Given a point $P\in M$, for all preimages of $p$ in $\tM$, consider the subset $P_j(i)$ such that $B(P_j(i),R) \subset A_i$. Denote $E_i$ a maximal set of $P_j(i)$ such that
\[
B(P_{j_1}(i),R)\cap B(P_{j_2}(i),R)=\emptyset
\]
for $j_1\neq j_2$. Following the argument of line 7 to line 23 on Page 152 of \cite{Liu2011short}, we can prove that there exists at least one $P_j(i)\in E_i$ such that
\begin{equation}\label{hgs:n-1}
\frac{\SL_r (B(P_j(i),R))}{\rvol(B(P_j(i),R))}\gs n-1-\Psi(\frac1i).
\end{equation}
Since
\begin{equation}\label{hls:n-1}
\Delta^{ac} r \ls (n-1)\frac{\cosh r}{\sinh r}=n-1+\Psi(r) \text{ and }\Delta^s r\ls 0.
\end{equation}
By combining (\ref{hgs:n-1}) with (\ref{hls:n-1}), we have
\begin{equation}\label{heq:n-1}
|\Delta^{ac} r-(n-1)| \ls \Psi(\frac1i) \text{ for } a.e. x \text{ and } -\Psi(\frac1i) \rvol(B(P_j(i),R))\ls \Delta^s r  (B(P_j(i),R))\ls 0.
\end{equation}

Fix $P_0\in \pi^{-1}(p)$, then there is an isometry $\Phi_i: B(P_0,R) \to B(P_j(i),R)$. Consider the function $u_i(x)=r(x)-|O P_j(i)|$ defined on $B(P_j(i),R)$. Let $v_i=u_i\circ \Phi_i: B(P_0,R) \to \bR$. $v_i$ is uniformly bounded and 1-Lipschitz, then there exists a subsequence (also denoted by $v_i$ for simplicity) uniformly converging to some $u_R$. Since $v_i$ is uniformly bounded in $W^{1,2}(B(P_0,R))$, by Lemma \ref{lem:Green}, for any $\varphi \in \rLip_c(B(P_0,R))$, we have
\begin{equation}\label{heq:LuR}
\begin{array}{ll}
\SL_{u_R}(\varphi)&=\lim_{i\to \infty} \SL_{v_i}(\varphi)\\
&=\lim_{i\to \infty} \SL_{u_i}(\varphi\circ \Phi_i^{-1})\\
\end{array}
\end{equation}
Denote $\varphi_i=\varphi\circ \Phi_i^{-1}$, by combing (\ref{heq:LuR}) with (\ref{heq:n-1}), we have
\begin{equation}
\begin{array}{ll}
\SL_{u_R}(\varphi)&=\lim_{i\to \infty} [\int_{B(P_j(i),R)} \varphi_i \Delta^{ac} r(x) \rd \rvol +\int_{B(P_j(i),R)} \varphi_i \rd \Delta^s r]\\
&=(n-1) \int_{B(P_j(i),R)} \varphi_i \rd \rvol\\
&=(n-1) \int_{B(P_0,R)} \varphi \rd \rvol.
\end{array}
\end{equation}
Since $v_i$ are $1+\Psi(\frac1i)$-concave, we know that $u_R$ is 1-concave. We claim that $|\nabla_x u_R|=1$ for a.e. $x\in B(P_0,R)$. In fact, Denote
\[
B_R'=\{x\in B(P_0,R)| v_i \text{ and }u_R \text{ are differentiable at } x \text{ for all } i\}.
\]
By Rademacher's theorem, it has full measure. For $x\in B'_R$, let $x_i=\Phi_i(x)$, choose geodesic connecting $x_i$ to $O$, let $\alpha_i=\Phi_i^{-1} \gamma_i$, then
\[
v_i(\alpha_i(t))-v_i(x)=u_i(\gamma_i(t))-u_i(x_i)=-t.
\]
Suppose that $\alpha_i$ subconverge to a geodesic $\alpha_R$, then
\[
u_R(\alpha_R(t))-u_R(x)=-t.
\]
This means that $|\nabla u_R(x)|=1$. So we prove the claim.

Suppose $u_R$ subconverge to some function $u:\tM \to R$. Repeat the above argument, we can prove that $u$ is 1-concave, $|\nabla_x u|=1$ for a.e. $x\in \tM$ and $\SL_u=(n-1)\cdot \rvol$. By Theorem \ref{thm:general}, we know that $M=\bR \times_{e^t} N$, where $N$ is an $n-1$ dimensional Alexandrov space with non-negative curvature. Then following the argument of the proof of Lemma 4.4 of Chen-Rong-Xu's paper \cite{CRX2016Quati}, $\tM$ is isometric to $\mathbb{H}^n$. For Reader's convenience, we list their argument below. Assume $\tM \ni \tilde{p}=(0,y)$ is a regular point, thus $\lim_{t\to \infty} (e^t N, y)=(\bR^{k-1},0)$. Via reparametrization of $s'=s-t$,
\begin{equation}
\begin{array}{ll}
\lim_{t\to \infty}(\bR \times_{e^s} N, (t,y))&=\lim_{t\to \infty} (\bR\times_{e^{s'}}e^t N, (0,y))\\
&=(\bR\times_{e^s} \bR^{k-1}, o)\\
&=(\bH^k,o).
\end{array}
\end{equation}
Since $M=\tM/ \pi_1(M)$ is compact, for any t, there is $\gamma_t\in \pi_1(M)$ such that 
\begin{equation}
|\gamma_t(\tilde{p}),(t,y)| \ls diam M \ls d.
\end{equation}
Then 
\begin{equation}
(\tM, \tilde{p})=\lim_{t\to\infty} (\tM, \gamma_t(\tilde{p}))=(\bH^k,o).
\end{equation}
\end{proof}

\end{document}